\documentclass[proc]{edpsmath}
\usepackage{dsfont}
\usepackage{color}
\usepackage{graphicx}
\usepackage{amsmath,amssymb,amsthm,textcomp}

\DeclareMathOperator\e{e}
\DeclareMathOperator\inter{int}

\begin{document}

\title{Numerical methods in the context of compartmental models in epidemiology}
\author{Peter Kratz}\address{Aix-Marseille Universit\'e, CNRS, Centrale Marseille, I2M, UMR 7373 13453
Marseille, France; \email{kratz@mathematik.hu-berlin.de; etienne.pardoux@univ-amu.fr; skbrice@yahoo.fr}}
\author{Etienne Pardoux}\sameaddress{1}
\author{Brice Samegni Kepgnou}\sameaddress{1}

%
%
\begin{abstract}
We consider compartmental models in epidemiology. For the study of the divergence of the stochastic model from its corresponding deterministic limit (i.e., the solution of an ODE) for long time horizon, a large deviations principle suggests a thorough numerical analysis of the two models. The aim of this paper is to present three such motivated numerical works. We first compute the solution of the ODE model by means of a non-standard finite difference scheme. Next we solve a constraint optimization problem via discrete-time dynamic programming: this enables us to compute the leading term in the large deviations principle of the time of extinction of a given disease. Finally, we apply the $\tau$-leaping algorithm to the stochastic model in order to simulate its solution efficiently. We illustrate these numerical methods by applying them to two examples.
\end{abstract}
\begin{resume}
On consid\`ere des mod\`eles comportementaux en \'epid\'emiologie. Afin d'\'etudier l'\'ecart en temps long entre le mod\`ele stochastique et sa limite loi des grands nombres (qui est la solution d'une EDO),  on se base sur un principe des grandes
d\'aviations, qui nous conduit \`a mener une \'etude num\'erique des deux mod\`eles, sur trois aspects diff\'erents. Tout d'abord, nous calculons une solution approch\'ee de l'EDO \`a l'aide d'une m\'ethode num\'erique dite ``non--standard''. Ensuite une r\'esolvons
num\'eriquement un probl\`eme de contr\^ole sous contrainte, afin de calculer le terme principal des grandes d\'eviations du temps de sortie d'une situation end\'emique. Enfin nous mettons en oeuvre l'algorithme du ``$\tau$--leaping'' pour simuler efficacement la solution du syst\`eme stochastique. Nous illustrons ces simulations num\'eriques en les appliquant \`a deux exemples.
 \end{resume}
\maketitle
\section*{Introduction}

It is well-known that deterministic ODE models of population dynamics (such as the evolution of diseases) are not appropriate for small population sizes $N$ (i.e., $N<10^3,10^6$). Recently, it has been shown (see~\cite{Campillo2012}) that in some cases the ODE models can diverge from the corresponding stochastic models even for larger population sizes ($N>10^6$).

We consider the ``natural'' stochastic extensions for many compartmental ODE models in epidemiology: individual based Poisson driven models. The ODE models can be recovered from these models as the law of large numbers limit as the population size $N$ approaches infinity. Hence, the deterministic model can be considered an appropriate first approximation of the stochastic model for large $N$. While the dynamics of many of the underlying deterministic models are well-understood, the analysis of the corresponding stochastic models is often difficult.

We are interested in the long-term behavior of the epidemic process. For many deterministic models, the disease will finally either die out or become endemic, depending on the parameter choice and/or the initial sizes of the compartments, i.e., the solution of the ODE converges to an asymptotically stable equilibrium. The theory of Large Deviations (also called the large deviations principle (LDP)) allows us to quantify the time taken by the stochastic model to diverge significantly from its deterministic limit. In particular, one can compute asymptotically, for large $N$, the expected time which is needed to transfer the epidemic process from the domain of attraction of one equilibrium to the one of another equilibrium. For the simplest models, this helps, e.g., to answer the question of when an epidemic becomes extinct although it should be endemic according to the deterministic approximation. Particularly interesting cases are epidemiological models where the ODE has several stable equilibria (usually a disease-free equilibrium and an endemic equilibrium as, e.g., in~\cite{Kribs2000,Chitnis2006}).

The remainder of this article is structured as follows. In section~\ref{SecModel}, we formulate stochastic and deterministic compartmental models and introduce two specific models which we analyze throughout this article by means of numerical methods (section~\ref{SubSecModel}). We furthermore describe how the time of exit from the domain of attraction of a stable equilibrium can be computed (section~\ref{SubSecLDP}). This motivates three different numerical projects which we outline in section~\ref{SubSecQuestions}. We address these questions in sections~\ref{SecODE}~-~\ref{SecSim} as follows. First we apply a non-standard finite difference scheme due to~\cite{Anguelov2014} in order to compute the solution of the deterministic model numerically. Second,  in section~\ref{SecVBar}, by solving a control problem via dynamic programming, we compute the leading coefficient in the large deviations analysis of the time of exit
 $\tau^N$ by the solution of the SDE from the domain of attraction of a locally stable equilibrium of the ODE. Finally, we simulate the stochastic process using the so-called \emph{$\tau$-leaping} algorithm.

\section{Compartmental models in epidemiology} \label{SecModel}

\subsection{Stochastic and deterministic models} \label{SubSecModel}

We consider diseases in large populations of (initially) $N$ individuals which are split up into different groups (or compartments) according to the disease status of the individuals. Such groups are, e.g., the group of individuals susceptible to the disease, the group of infectious individuals and the group of immune individuals. In order to better compare the population dynamics for different values of $N$, we normalize the compartment sizes by dividing by $N$; instead of considering the number of individuals in the different groups, we hence rather consider their ``proportions'' with respect to the total size of the population (note that for models with non-constant population size, these are only initially truly proportions).

For a closed set $A\subset \mathds{R}^d$, we define the $d$-dimensional jump process $Z^N$ via a set of $k$ jump directions $h_j\in \mathds{Z}^d$ and respective Poisson rates $\beta_j(z)$, $j=1,\dots,k$; $Z^N_i(t)$ denotes the ``proportion'' of individuals in compartment $i=1,\dots,d$ at time $t\geq 0$:
\begin{align}\label{EqPoisson}
Z^N(t):=Z^{N,x}(t)  &:= x +\frac{1}{N} \sum_{j=1}^k h_j P_j\Big(\int_0^t N \beta_j(Z^N(s)) ds \Big) \\
&= x +\int_0^t \sum_{j=1}^k h_j \beta_j(Z^{N}(s)) ds+ \frac{1}{N} \sum_{j=1}^k h_j M_j\Big(\int_0^t N \beta_j(Z^N(s)) ds \Big); \notag
\end{align}
here $x\in A$,  $(P_j)_{j=1,\dots,k}$ are i.i.d.~standard Poisson processes and $M_j(t)=P_j(t)-t$, $j=1,\dots,k$, are the associated martingales. The jump directions $h_j$ and their respective rates are chosen in such a way that $Z^N(t)\in A$ almost surely for all $t$.

The dynamics of the deterministic model corresponding to the stochastic process in~\eqref{EqPoisson} are given by the solution of the following ODE:
\begin{equation}\label{EqODE}
Z(t):=Z^{x}(t)  = x +\int_0^t \sum_{j=1}^k h_j \beta_j(Z(s)) ds.
\end{equation}

The two models are connected via a Law of Large Numbers (LLN) (see~\cite{Kurtz1978}): if the rates $\beta_j$ are Lipschitz continuous, we have
\begin{equation} \label{EqLLN}
\lim_{N \rightarrow \infty} Z^{N,x}(t) = Z^x(t) \quad \text{a.s. uniformly on compact time intervals } [0,T].
\end{equation}
Note also that in the case of Lipschitz continuous rates, equation~\eqref{EqODE} admits a unique solution on $[0,T]$.

\begin{xmpl} \label{ExampleSIS}
\begin{bfseries} SIS model.\end{bfseries}
We first consider a simple SIS model without demography ($S(t)$ being the number of susceptible individuals and $I(t)=N-S(t)$ the number of infectious individuals at time $t$) . For $\beta>0$, we assume that the rate of infections is $\beta S(t) I(t)/N$.\footnote{The reasoning behind this is the following. Each individual (in particular each infected individual) meets other individuals at rate $\alpha$. The probability that the encounter is with a susceptible is $S(t)/N$. Denote by $p$ the probability that such an event yields a new infection. Hence the total rate of new infections is $\beta S(t) I(t)/N$, if $\beta=p\alpha$.} In a similar way, we assume for $\gamma>0$ that the rate of recoveries is $\gamma I(t)$. We ignore immunity due to earlier infection and assume that an infectious individual becomes susceptible after recovery. As population size is constant, we can reduce the dimension of the model by solely considering the proportion of infectious at time $t$,
that is $d=1$ and $Z_t=I_t$. Using the notation of equations~\eqref{EqPoisson} and~\eqref{EqODE}, we have
\[
A=[0,1], \quad h_1=1, \quad  h_2=-1, \quad \beta_1(z)= \beta z (1-z), \quad \beta_2(z)=\gamma z.
\]
It is easy to see that the ODE~\eqref{EqODE} has a disease free equilibrium $\bar x=0$. This equilibrium is asymptotically stable if $R_0=\beta/\gamma<1$.\footnote{$R_0$ is the so-called \emph{basic reproduction number}. It denotes the average number of secondary cases infected by one primary case during its infectious period, at the start of the epidemic.}
If $R_0>1$, $\bar x$ is unstable and there exists a second, endemic equilibrium $x^*=1-\gamma/\beta$ which is asymptotically stable. While in the deterministic model the proportion of infectious individuals converges to the endemic equilibrium $x^*$, the disease goes ultimately extinct in the stochastic model.
\end{xmpl}

\begin{xmpl} \label{ExampleSIV}
\begin{bfseries}A model with vaccination.\end{bfseries} We consider a (deterministic) model with vaccination and demography described in~\cite{Kribs2000} and its stochastic counterpart. We illustrate the different transitions in Figure~1; here, $S$ and $I$ denote the number of susceptible respectively infectious individuals as before and $V$ denotes the number of vaccinated individuals. We assume that individuals are vaccinated at a certain rate but can lose their protection again; the vaccine is not perfect but decreases the rate of infection by a factor $\sigma\in [0,1]$: if $\sigma=1$, the vaccine is useless, if $\sigma=0$ it protects the individuals perfectly. Furthermore, we assume that individuals are born susceptible and die (at the same rate) independently of their disease status. This ensures that the population size remains constant in the deterministic model. For simplicity, we ensure constant population size by synchronizing births and deaths in the stochastic model. Hence, we can again reduce the dimension of the models; the first coordinate of the process denotes the proportion of infective individuals, the second coordinate denotes the proportion of vaccinated individuals. Using the notations of equations~\eqref{EqPoisson} and~\eqref{EqODE}, we obtain
\begin{align*}
&A = \{z \in \mathds{R}^2_+ | 0 \leq z_1+z_2 \leq 1 \},  \\
&h_1= (1,0)^\top, \quad \beta_1(z)=\beta z_1 (1-z_1-z_2), \quad
h_2= (1,-1)^\top, \quad \beta_2(z)=\sigma\beta z_1 z_2, \\
&h_3= (-1,0)^\top, \quad \beta_3(z)= \gamma z_1, \quad
h_4= (0,1)^\top, \quad \beta_4(z)=\theta z_2, \quad
h_5= (0,-1)^\top, \quad \beta_5(z)=\eta (1-z_1- z_2), \\
&h_6= (-1,0)^\top, \quad \beta_6(z)= \mu z_1, \quad
h_7= (0,-1)^\top, \quad \beta_7(z)= \mu z_2.
\end{align*}
The basic reproduction number of this model is given by
\begin{equation*}
  R_{0}=\frac{\beta}{\mu+\gamma}\frac{\mu+\theta+\sigma\eta}{\mu+\theta+\eta}, \quad
  \tilde{R}_{0}:=\frac{\beta}{\mu+\gamma};
\end{equation*}
here, $\tilde R_0$ denotes the basic reproduction number without vaccination, i.e., $V(0)=0$ and $\phi=0$ (see~\cite{Kribs2000}). There exists a disease-free equilibrium $\bar x$ ($\bar x_1=0$) which is asymptotically stable for $R_{0}<1$  and unstable
for $\tilde{R}_{0}>1$. For $R_0 < 1 < \tilde R_0$ (and additional assumptions on the parameters, see~\cite{Kribs2000} again), two endemic equilibria $x^*$ and $\tilde x$ ($x^*_1, \tilde x_1>0$) exist: one (say $x^*$) is locally asymptotically stable and one is unstable; the disease-free equilibrium is locally asymptotically stable.
\end{xmpl}

\begin{figure} 
\begin{picture} (290,140)
\linethickness{0.4mm}

\put(135,0){\line(1,0){40}}
\put(135,0){\line(0,1){40}}
\put(175,0){\line(0,1){40}}
\put(135,40){\line(1,0){40}}

\put(60,90){\line(1,0){40}}
 \put(60,90){\line(0,1){40}}
\put(100,90){\line(0,1){40}}
\put(60,130){\line(1,0){40}}

\put(200,90){\line(1,0){40}}
\put(200,90){\line(0,1){40}}
\put(240,90){\line(0,1){40}}
\put(200,130){\line(1,0){40}}

\put(190,95){\vector(-1,0){80}}
\put(110,125){\vector(1,0){80}}
\put(70,85){\vector(3,-4){55}}
\put(180,10){\vector(3,4){55}}
\put(130,40){\vector(-3,4){30}}

\put(250,110){\vector(1,0){40}}
\put(45,95){\vector(-1,0){40}}
\put(5,125){\vector(1,0){40}}
\put(180,5){\vector(1,0){40}}

\put(72,100){\Huge{$S$}}
\put(215,100){\Huge{$I$}}
\put(144,10){\Huge{$V$}}

\put(120,129){\small{$\beta S I/N$}}
\put(120,115){\small{(infection rate)}}

\put(120,99){\small{$\gamma I$}}
\put(120,85){\small{(recovery rate)}}

\put(65,50){\small{$\eta S$}}
\put(20,40){\small{(vaccination rate})}

\put(115,65){\small{$\theta V$}}
\put(123,55){\small{(loss of vaccination)}}

\put(15,129){\small{$\mu N$}}
\put(5,115){\small{(birth rate)}}

\put(15,99){\small{$\mu S$}}
\put(5,85){\small{(death rate)}}

\put(265,114){\small{$\mu I$}}
\put(220,50){\small{$\sigma \beta V I/N \,\, (\sigma \in [0,1])$}}

\put(215,40){\small{(infection of vaccinated)}}
\put(200,10){\small{$\mu V$}}

\end{picture}
\caption{Schematic representation of the disease transmission for the model with vaccination by~\cite{Kribs2000}.}
\end{figure}
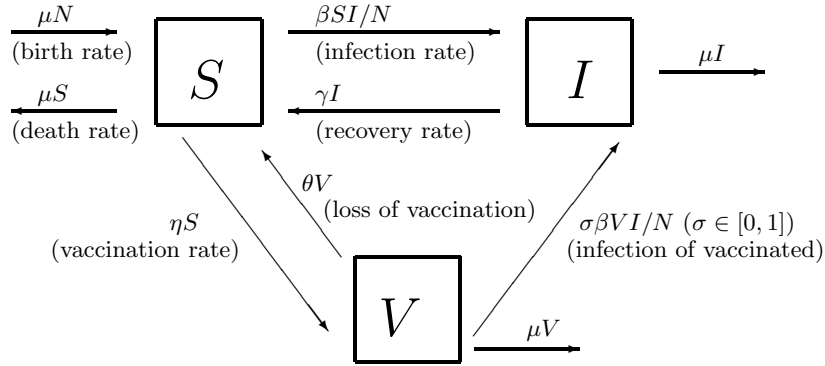

\subsection{Large deviations} \label{SubSecLDP}

In order to quantify the deviation of the stochastic model from the deterministic model, an LDP is desirable. For the models we have in mind (cf.~Examples~\ref{ExampleSIS} and~\ref{ExampleSIV}), some of the rates $\beta_j:A \rightarrow \mathds{R}_+$ tend to zero as $x$ approaches the boundary of $A$ (this ensures that the process does not leave the domain $A$). For this reason, general results for LDPs for Markov jump processes (see \cite{Shwartz1995,Dupuis1997,Feng2006}) do not apply. \cite{Shwartz2005} provides a large deviations principle which allows for such diminishing rates; however, their assumptions do often not apply to the more complicated models in epidemiology.\footnote{The assumptions in~\cite{Shwartz2005} apply to the SIS model (Example~\ref{ExampleSIS}) but not to the model with vaccination in Example~\ref{ExampleSIV}.} \cite{Kratz2014} provides an LDP on the space $D([0,T];A)$\footnote{Here, $D([0,T];A)$ denotes the space of c\`{a}dl\`{a}g functions
on $[0,T]$ with domain $A$ equipped with an appropriate metric such that the resulting space is Polish (see, e.g., \cite{Billingsley1999}).} for a large class of epidemiological models with compact domain $A$ by generalyzing the results from~\cite{Shwartz2005}; the rate function $I_{T,x}$ of the LDP is given as follows:\footnote{See, e.g., \cite{Dembo2009} for an introduction to the theory of large deviations and a definition of the rate function.}
\begin{align*}
I_{T,x}(\phi)&:= \begin{cases}
\int_0^T L(\phi(t),\phi'(t))  dt& \text{ if } \phi \text{ is absolutely continuous and } \phi(0)=x \\
\infty & \text{ else,}
\end{cases}
\end{align*}
where $L$ denotes the Legendre-Fenchel transform:
\begin{equation} \label{EqLegendre}
L(x,y):=\sup_{p \in \mathds{R}^d} \ell (p,x,y)
\end{equation}
for
\[
\ell(p,x,y) = \langle p,y \rangle - \sum_{j=1}^k \beta_j(x) ( \e^{\langle p, h_j \rangle} -1).
\]
By~\cite{Shwartz2005}, we have that
\begin{equation} \label{EqODE=0}
I_{T,x}(\phi)=0 \quad \text{if and only if $\phi$ solves the ODE~\eqref{EqODE}}.
\end{equation}
We interpret $I_{T,x}(\phi)$ as the \emph{action} required to follow a trajectory $\phi$, or in other words the action required to deviate from the solution of the ODE.

We are particularly interested in the Freidlin-Wentzell theory (see, e.g.,~\cite{Freidlin1998,Dembo2009}). First, we want to compute the \emph{time} $\tau^{N,x}$ at which the process $Z^{N,x}$ leaves the domain of attraction of an asymptotically stable equilibrium $O$, i.e.,
\[
\tau^{N,x}:=\inf \{ t >0| Z^{N,x}(t) \in A\setminus O\}.
\]
In the case of the SIS model (and $R_0>1$) we consider the set $O=(0,1]$; for the model with vaccination (and appropriate parameter choice, in particular $R_0 < 1 < \tilde R_0$), we consider the set
\[
O=\{z \in A| \lim_{t \rightarrow \infty} Z^{z}(t) = x^* \}. \footnotemark
\]
\footnotetext{We want to remark here that we cannot consider the exit from the domain of attraction of the disease-free equilibrium $\bar x$, i.e.,
\[
O=\{z \in A| \lim_{t \rightarrow \infty} Z^{z}(t) = \bar x \}.
\]
The reason for this is that the absorbing set $\tilde A=\{z\in A| z_1=0\}$ (where $z_1$ is the first coordinate of $z$) cannot be left both by the solution of the ODE and by the solution of the SDE.}
In other words, we are interested in time of extinction of the disease\footnote{Respectively in the time when the process is no longer attracted by the endemic equilibrium and hence the probability of extinction is significantly increased, cf.~the LLN, equation~\eqref{EqLLN}.} which should eventually become endemic according to the deterministic model.

If $A$ is compact, it can be deduced from the LDP that for all $\delta>0$,
\begin{equation} \label{EqExitTime}
\lim_{N\rightarrow \infty} \mathbb{P}\big[ \e^{N(\bar V - \delta)} < \tau^{N,x} < \e^{N(\bar V + \delta)} \big]=1,
\end{equation}
where $\bar V$ is given by the solution of the following optimization problem:
\begin{equation} \label{EqVBar}
\bar V:= \inf_{z \in A\setminus O} V(x^*,z)
\end{equation}
\[
\text{for} \quad V(x^*,z):=\inf_{T>0, \phi \in D([0,T];A): \phi(0)=x^*, \phi(T)=z} I_{T,x^*}(\phi)
\]
(see~\cite{Kratz2014}).

Second, we are interested in the \emph{place}, where $Z^N$ leaves the domain $O$. Although existing Freidlin-Wentzell theory cannot be applied to the epidemiological models in Examples~\ref{ExampleSIS} and~\ref{ExampleSIV}, we strongly believe that the following holds (see~\cite{Freidlin1998,Dembo2009} again for corresponding results for other processes): for any compact set $C \subset \partial O\cap (A \setminus O)$ with $\bar V < \inf_{z \in C} V(x^*,z)$, $x \in O$
\[
\lim_{N\rightarrow \infty} \mathbb{P}\big[Z^{N,x}(\tau^{N,x})\in C\big]=0.
\]
In particular, if there exits a $z^*\in \partial O\cap (A \setminus O)$ with $V(x^*,z^*)<V(x^*,z)$ for all other $z\in \partial O\cap (A \setminus O)$, then for all $\delta>0$,
\begin{equation} \label{EqExitPlace}
\lim_{N\rightarrow \infty} \mathbb{P}\big[|Z^{N,x}(\tau^{N,x})-z^*|>\delta]=0.
\end{equation}

\subsection{Research questions} \label{SubSecQuestions}

In order to better understand the differences between deterministic and stochastic compartmental models in epidemiology, the theoretical results above suggest a thorough numerical analysis of the specific models of Examples~\ref{ExampleSIS} and~\ref{ExampleSIV}. In the remainder of this article, we address the following questions.

\begin{enumerate}
\item[1.]
In order to analyze the deviation of the stochastic model from its Law of Large Numbers limit, we first have to compute the solution of the deterministic models in a reliable way in section~\ref{SecODE}. In order to avoid instabilities of the solution (such as components becoming negative), we apply a \emph{non-standard finite difference scheme} developed in~\cite{Anguelov2014}.
\item[2.]
The computation of the time (and place) of exit from domain depends on the solution of the optimization problem~\eqref{EqVBar}. We use dynamic programming to solve the optimization problem numerically and hence deduce the results about the time of exit according to~\eqref{EqExitTime}   in section~\ref{SecVBar}.
\item[3.]
Finally, we simulate the process $Z^N$ in section~\ref{SecSim}. It is not difficult to simulate the process in an exact way by using Gillespie's stochastic simulation algorithm (SSA) (see, e.g., \cite{Gillespie1976, Gillespie1977}). However, the simulation is rather slow  as soon as $N$ becomes large, and we hence adapt the \emph{$\tau$-leaping} algorithm (see, e.g.,~\cite{Gillespie2001}).
\end{enumerate}

\section{Computation of the solution of the ODE} \label{SecODE}
The large deviations theory will tell us about the deviations between the long--term behaviors of the solution  of the SDE and of its ODE LLN limit.
Hence we are interested in the long-term behavior of the deterministic epidemiological models described by equation~\eqref{EqODE}. Consequently we need a numerical scheme to compute the solution of the ODE~\eqref{EqODE} in a ``reliable'' way. By reliable, we mean in particular that the fixed points of the solution of the (discrete) numerical scheme are the same as the equilibria of the ODE (with the same local stability), which forces the long time bahavior of the numerical scheme to resemble that of the solution of the ODE. Explicit finite difference schemes for the solution of ODEs can have undesirable instabilities. In section~\ref{SubSecODESIS} below, we consider the SIS model from Example~\ref{ExampleSIS} and illustrate that the solution of such a scheme can be oscillating around zero (and in particular become negative) if the parameters (particularly the step-size $h$) are not chosen with great care. For this reason, we apply a non-standard finite difference (NSFD) method for models in epidemiology which respects the qualitative behavior of the solution of the ODE (see~\cite{Anguelov2014}). This NSFD method follows two basic rules for NSFD methods (for a detailed description of the NSFD method see, e.g., \cite{Mickens2000, Mickens2005}). First, quadratic terms of the ODE are approximated in a non-local way (by mixed explicit/implicit terms). Second, non-trivial denominator functions are used for the discrete derivatives.

The remainder of this section is structured as follows. In section~\ref{SubSecODENSFD}, we describe the specific NSFD method for epidemiological models from~\cite{Anguelov2014} and outline in which sense the qualitative long-term behavior of the ODE is respected. In section~\ref{SubSecODESIS}, we apply the method to the SIS model (Example~\ref{ExampleSIS}) and compare it to the corresponding explicit scheme. Finally, we apply the method to the model with vaccination (Example~\ref{ExampleSIV}) and compute the boundary separating the domains of attraction of the stable endemic equilibrium and the disease-free equilibrium.

\subsection{A NSFD method for epidemiological models} \label{SubSecODENSFD}

We consider deterministic models for which the ODE~\eqref{EqODE} can be rewritten in the following form
\begin{equation}\label{EqODEMetzler}
\frac{dZ}{dt} =A(Z)Z+f, \quad Z(0)= x,
\end{equation}
where $A(Z)=(a_{i,j}(Z))_{i,j=1,\dots,d}\in \mathbb{R}^{d \times d}$ is a non-linear Metzler matrix\footnote{A Metzler matrix is a matrix  $A=(a_{i,j})_{i,j=1\dots,d}$ which satisfies $a_{i,j}\geq 0$ for $i\neq j$.} and $f \in \mathbb{R}^d$; usually, $f$ represents the rates of birth and immigration into the different compartments (if they exist, cf.~Example~\ref{ExampleSIV}).

We follow~\cite{Anguelov2014} and discretize time by fixing the time step $h>0$; we write $t_{m}=mh$ for $m\in\mathbb{N}$ and approximate the solution $Z$ of the ODE~\eqref{EqODEMetzler} at the discrete time points $t_m$ by the solution $\tilde Z$ of the following implicit NSFD scheme (i.e., $\tilde{Z}^{m}\approx Z(t_{m})$):
\begin{equation} \notag
\frac{\tilde{Z}_{i}^{m+1}-\tilde{Z}_{i}^{m}}{\psi(h)} =\sum_{j=1}^{d}a_{i,j}(\tilde{Z}^{m})\tilde{Z}_{j}^{m+1}+f_{i} \quad (i=1,\dots,d), \quad \tilde{Z}^{0}= x;
\end{equation}
the denominator function $\psi$ is specified below and satisfies $\psi(h)=h+o(h^{2})$. We first note that this scheme can equivalently be written as
\begin{equation}
(I-\psi(h)A(\tilde{Z}^{m}))\tilde{Z}^{m+1} =\tilde{Z}^{m}+\psi(h)f, \quad
\tilde{Z}^{0}= x;
\label{EqNSFDMatrix}
\end{equation}
$(I-\psi(h)A(\tilde{Z}^{m}))$ is an $M$-matrix\footnote{An $M$-matrix is a matrix $A$ with $a_{i,j}\leq 0$ for $i\neq j$ which is invertible such that all entries of its inverse are non-negative.} and hence invertible (see~\cite{Anguelov2014}). For the definition of the denominator function, we assume that all equilibria of the ODE~\eqref{EqODEMetzler} are hyperbolic. With $\varphi:\mathbb{R}\rightarrow\mathbb{R}$ satisfying $\varphi(h)=h+o(h^{2})$ and $0<\varphi(z)<1$ for $z>0$ (e.g., $\varphi(z)=1-\exp(-z)$ or $\varphi(z)=z/(1+z^{2})$), $\psi$ is  defined as
\begin{equation} \label{EqFuncDem}
\psi(h):=\frac{\varphi(Qh)}{Q} \quad \text{for} \quad Q\geq\max\left\{\frac{|\lambda|^{2}}{2|Re(\lambda)|}\right\},
\end{equation}
where the maximum is taken over the eigenvalues $\lambda$ of the Jacobian of the right hand side of \eqref{EqODEMetzler} at the equilibrium points. In particular, we have $\psi(h)=h+o(h^{2})$ as desired.

The NSFD scheme~\eqref{EqNSFDMatrix} is \emph{elementary stable} in the sense that it has exactly the same fixed points as the continuous ODE~\eqref{EqODE} it approximates, with the same local stability for all values of $h$ (see~\cite{Anguelov2014}). In addition, it can be ensured that the components of the solution of the scheme remain positive by choosing $h$ and hence $\phi(h)$ small enough (see~\cite{Anguelov2014} again). For the examples we consider, the solution of the scheme always stays in the ``natural domain'' of the solution of the ODE~\eqref{EqODE}.

\subsection{The SIS model} \label{SubSecODESIS}

For the SIS model of Example~\ref{ExampleSIS}, the ODE~\eqref{EqODE} becomes
\begin{equation} \label{EqODESIS}
Z(t)=x+  \int_0^t \big(\beta Z(s) (1-Z(s)) - \gamma Z(s) \big) ds =x + \int_0^t \big( - \beta  Z(s)^2 + (\beta-\gamma) Z(s) \big) ds;
\end{equation}
recall that in this case $Z(t)$ denotes the number of infectious individuals at time $t$. This is a Riccati differential equation with constant coefficients which admits an explicit solution:
\begin{equation} \label{EqODESISExact}
Z(t)=
\begin{cases}
\frac{(\beta-\gamma)x \exp((\beta-\gamma) t)}{(\beta-\gamma)+ \beta x (\exp((\beta-\gamma) t)-1)} &\text{if } \beta\not=\gamma \\
\frac{x}{1-\beta x t} &\text{else.}
\end{cases}
\end{equation}
We nevertheless compute the solution of the ODE~\eqref{EqODESIS} numerically by a standard (explicit) difference scheme and by the NSFD scheme introduced in section~\ref{SubSecODENSFD} in order to illustrate possible instabilities of the explicit scheme.

Let us first consider the following explicit scheme:
\begin{equation} \label{EqODESISSFD}
\frac{\tilde Z^{m+1}-\tilde Z^{m}}{h} = \beta \tilde Z^{m}(1-\tilde Z^{m})-\gamma \tilde Z^{m}, \quad \tilde Z^m=x
\end{equation}
or equivalently
\[
\tilde Z^{m+1} = \tilde Z^{m}(1-\gamma h+\beta h)-h\beta(\tilde Z^{m})^{2}, \quad \tilde Z^m=x.
\]
We note that if, e.g., $\beta=1/h$ and $\gamma=2/h$, $\tilde Z^{m}<0$ for all $m\geq 1$ which is obviously undesirable. If, e.g., $\beta=1/h$ and $\gamma=4/h$, the solution is oscillating around zero with $\tilde Z^m<0$ for $m$ odd and $\tilde Z^m>0$ for $m$ even. We illustrate the solution for the latter case in the left picture of Figure~2.

Let us now apply the NSFD scheme to the model. We represent the ODE in the form~\eqref{EqODEMetzler} by setting $A(Z)=-\beta Z +\beta - \gamma$ and define the denominator function $\psi$ according to equation~\eqref{EqFuncDem} (for $\varphi(z)=1-\exp(-z)$). The solution does not suffer from the same instabilities as the solution of the explicit scheme as illustrated in the middle picture of Figure~2. We compare the solutions of the two schemes with the exact solution of the ODE as given in equation~\eqref{EqODESISExact} (right picture of Figure~2).
We want to remark here that the NSFD scheme is of course by no means the only possible scheme which avoids instabilities in this simple case. However,\cite{Anguelov2014} provides us with a theoretical basis which ensures that instabilities do not occur even for much more complicated models.

\begin{figure} 
\centering \includegraphics[height=0.225\hsize]{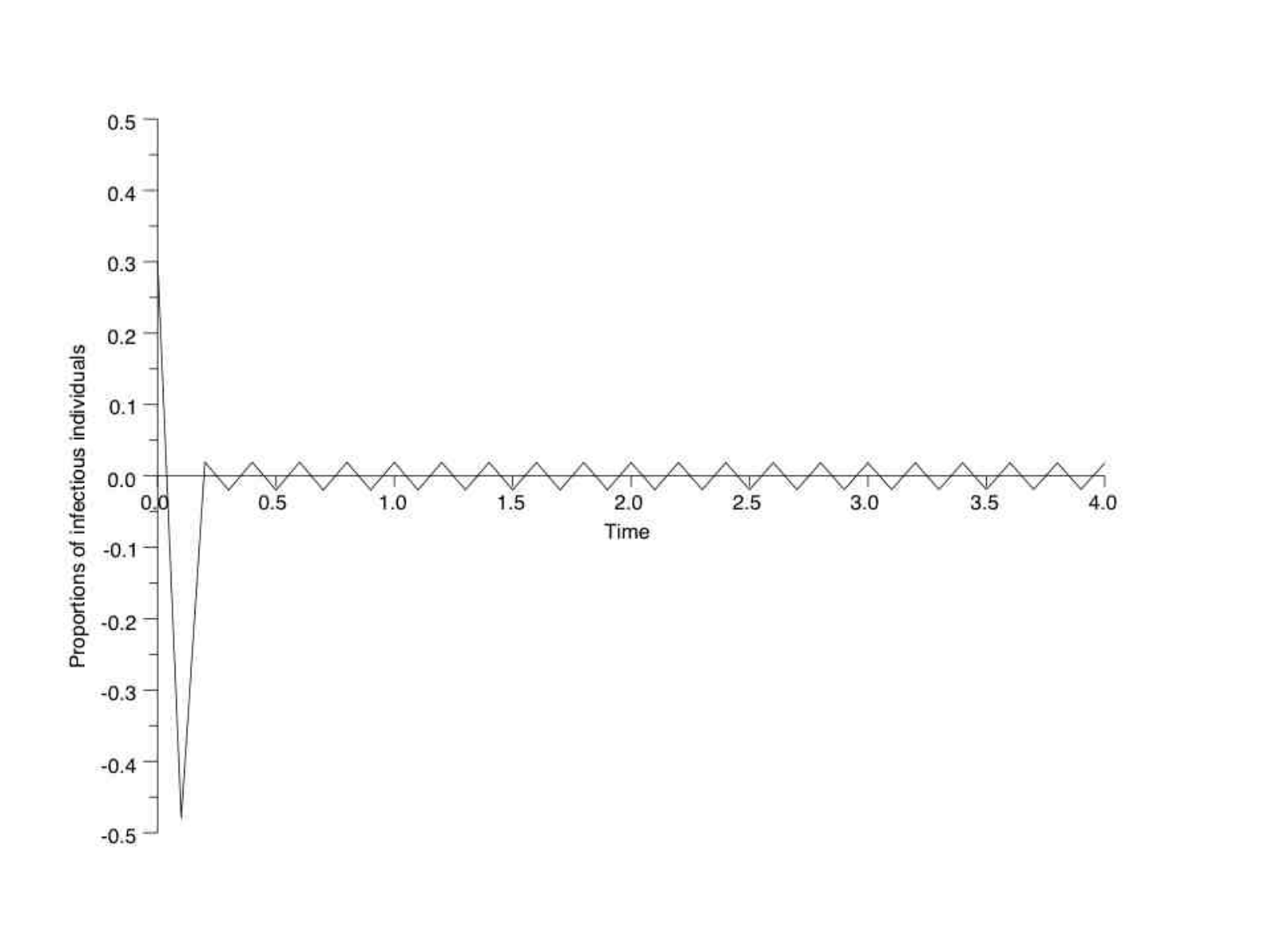}
\qquad
\centering \includegraphics[height=0.225\hsize]{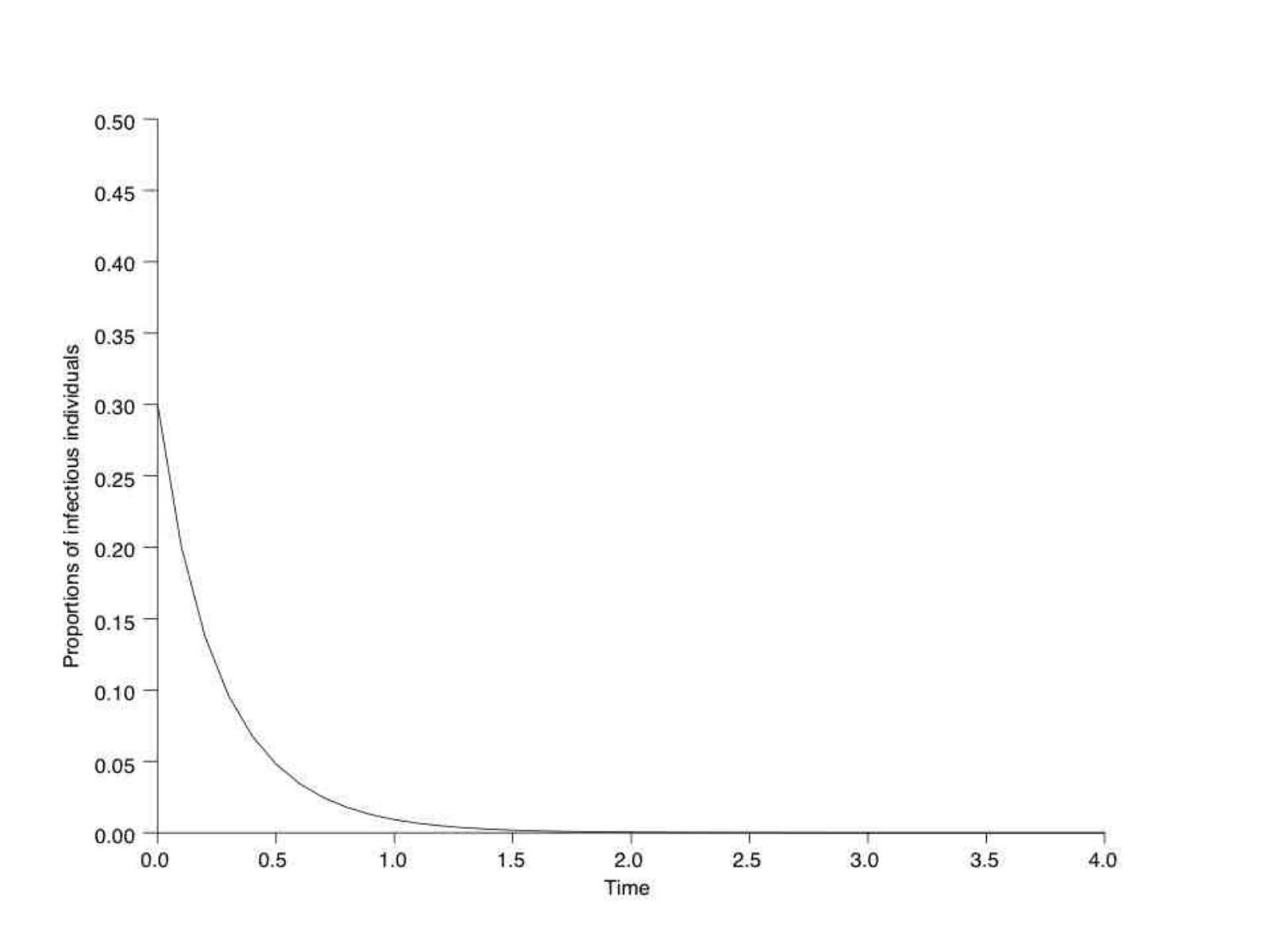}
\qquad
\includegraphics[height=0.225\hsize]{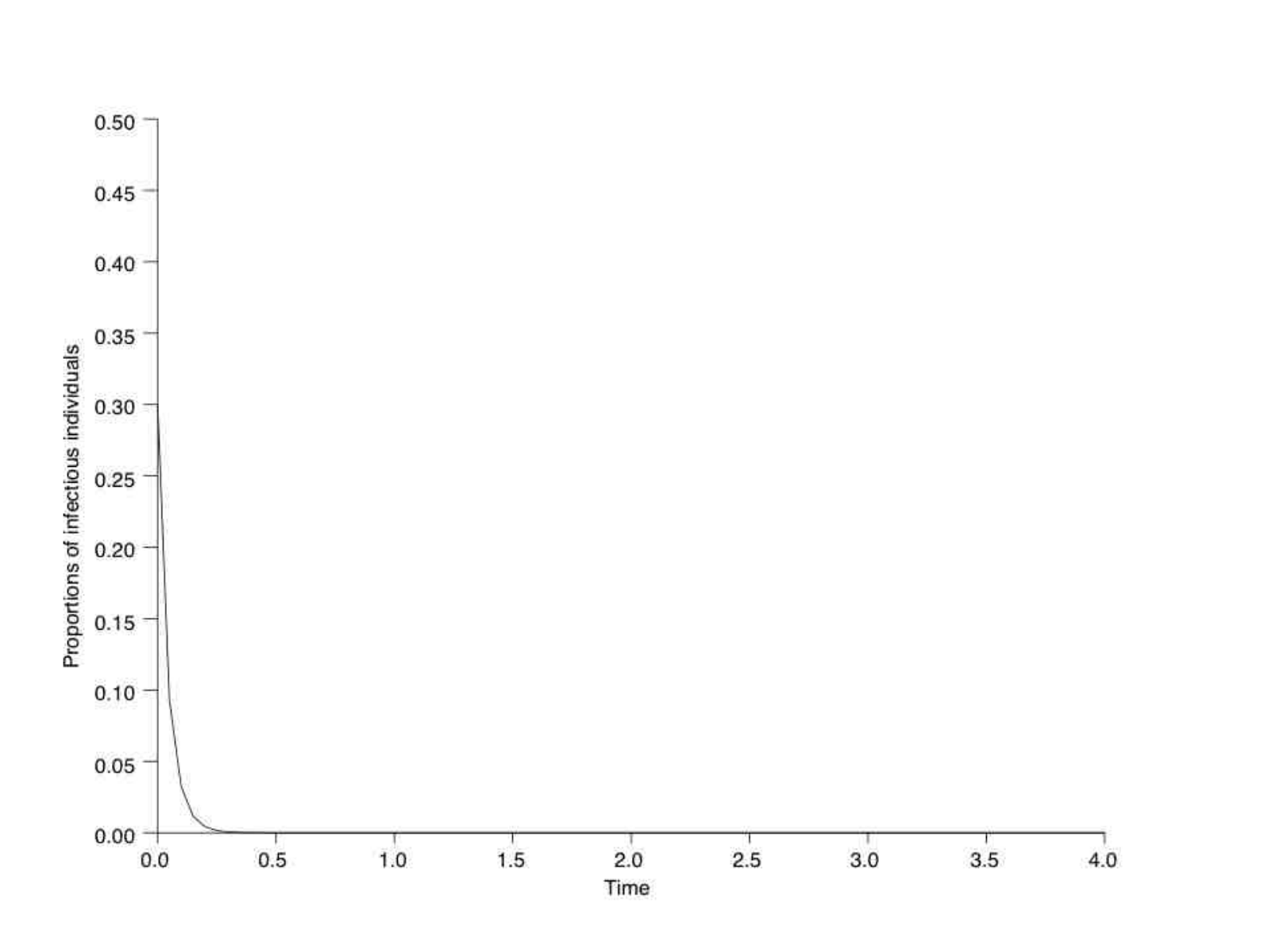}
\caption{The left picture denotes the solution of the explicit finite difference scheme~\eqref{EqODESISSFD}. The middle picture denotes the solution of the NSFD scheme. The right picture denotes the exact solution of the ODE according to equation~\eqref{EqODESISExact}. $Z(0)=0.3$, $T=4$, $h=0.1$, $\beta=4/h=40$, $\gamma=2/h=20$. We choose here a large step size $h$ on purpose, to show how the two schemes behave for large step size. Of course, both improve if we reduce $h$.}
\end{figure}

\subsection{The model with vaccination} \label{SubSecODESIV}

We now consider the model with vaccination from Example~\ref{ExampleSIV}. We denote the proportions of susceptible, vaccinated and infectious individuals at time $t$ by $Z_1(t)$, $Z_2(t)$ and $Z_3(t)=1-Z_1(t)-Z_2(t)$, respectively. Note that we cannot reduce the dimension from three to two if we want to represent the ODE in the form~\eqref{EqODEMetzler} (with appropriate Metzler matrix $A$). We set
\begin{equation*}
A(Z)=\left(
\begin{array}{ccc}
  -\beta Z_3-\mu-\eta & \theta & \gamma \\
  \eta & -\sigma\beta Z_3 & 0 \\
  \beta Z_3 & \sigma\beta Z_3 & -\mu-\gamma \\
\end{array}
\right)\quad \text{and} \quad f=(\mu,0,0)^\top
\end{equation*}
and choose the denominator function $\psi$ according to equation~\eqref{EqFuncDem} for $\varphi(z)=1-\exp(-z)$. We readily observe that the relevant assumptions for the NSFD scheme are satisfied and apply the scheme for time step $h=0.1$ to different initial values $Z(0)$ and a set of parameters (cf.~Figure~3) such that there exists a stable disease-free equilibrium ($\bar x =(0.14,0.86,0)^\top$), a stable endemic equilibrium ($x^*=(0.24,0.45,0.31)^\top$) and an unstable endemic equilibrium ($\tilde x=(0.23,0.59,0.18)^\top$; cf.~the corresponding discussion in Example~\ref{ExampleSIV}). We illustrate the evolution of the proportions of vaccinated and infectious individuals for different initial values in Figure~3. While the disease  eventually dies out in the left picture, it becomes endemic in the right picture.

\begin{figure} 
\centering \includegraphics[height=0.35\hsize]{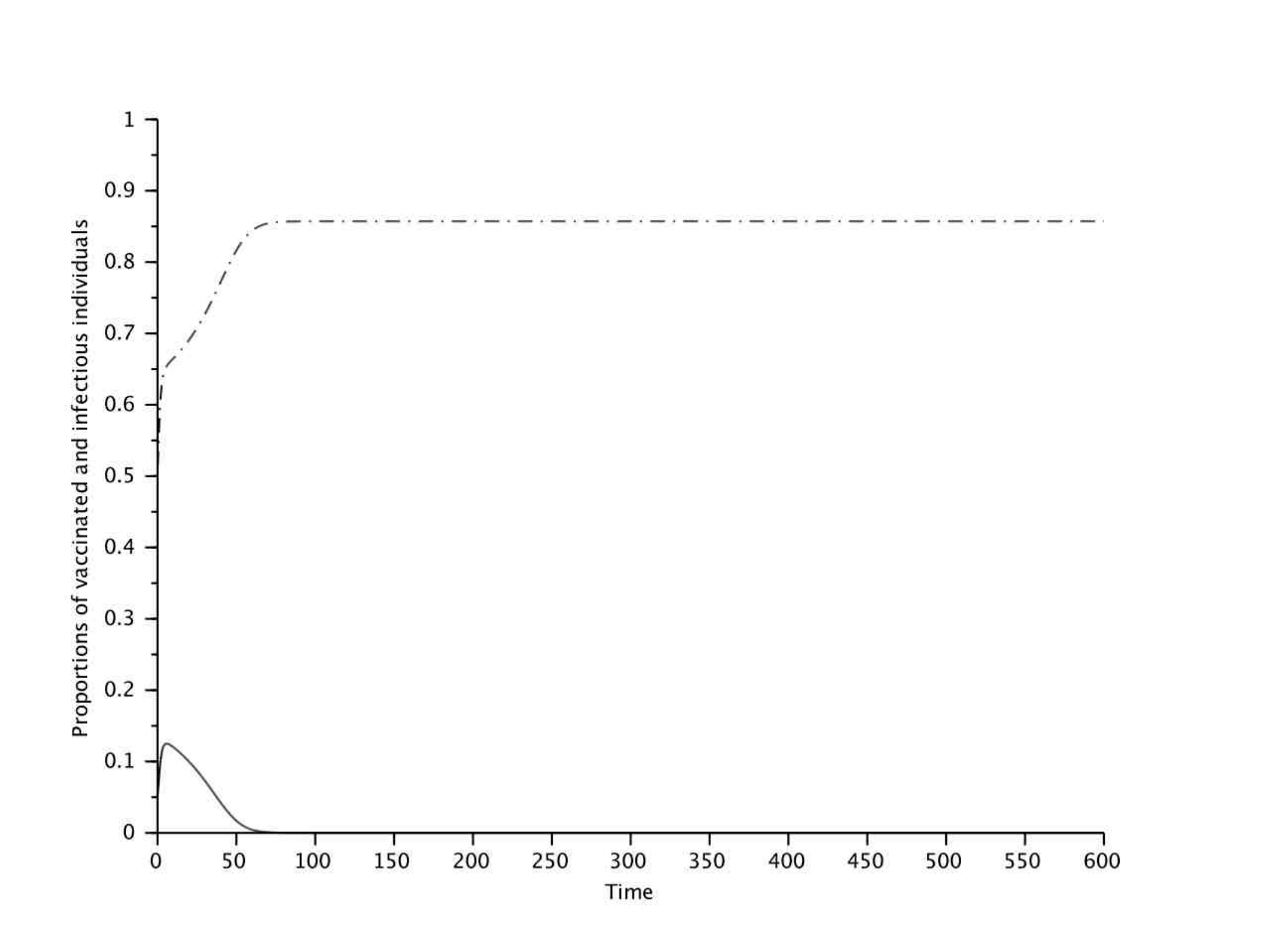}
\qquad
\centering \includegraphics[height=0.35\hsize]{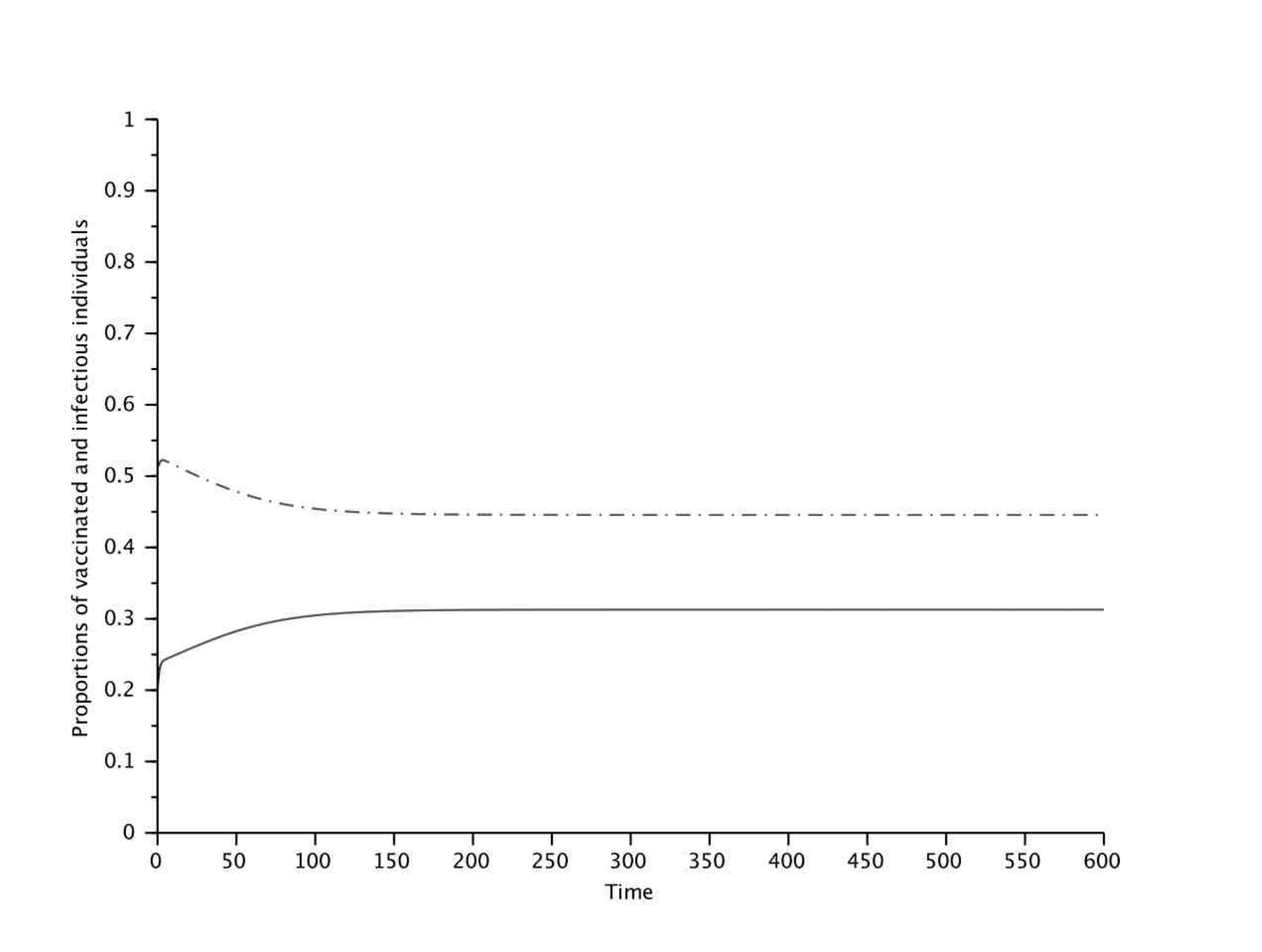}
\caption{Evolution of the proportions of vaccinated individuals (dashed lines) and infectious individuals (solid lines) for initial values $Z_2(0)$=0.5, $Z_3(0)=0.05$ (left picture) respectively $Z_2(0)$=0.5, $Z_3(0)=0.2$ (right picture). $T=600$, $\beta=3.6$, $\gamma=1$, $\eta=0.3$, $\theta=0.02$, $\mu=0.03$, $\sigma=0.1$.}
\end{figure}

The NSFD scheme also allows us to approximately compute the \emph{characteristic boundary}, i.e., the boundary separating the domains of attractions of the two stable equilibria $\bar x$ and $x^*$ (in other words: the set of points for which the solution of the ODE converges to the unstable equilibrium $\tilde x$). We illustrate the characteristic boundary in Figure~4.

\begin{figure} 
\includegraphics[height=0.35\hsize]{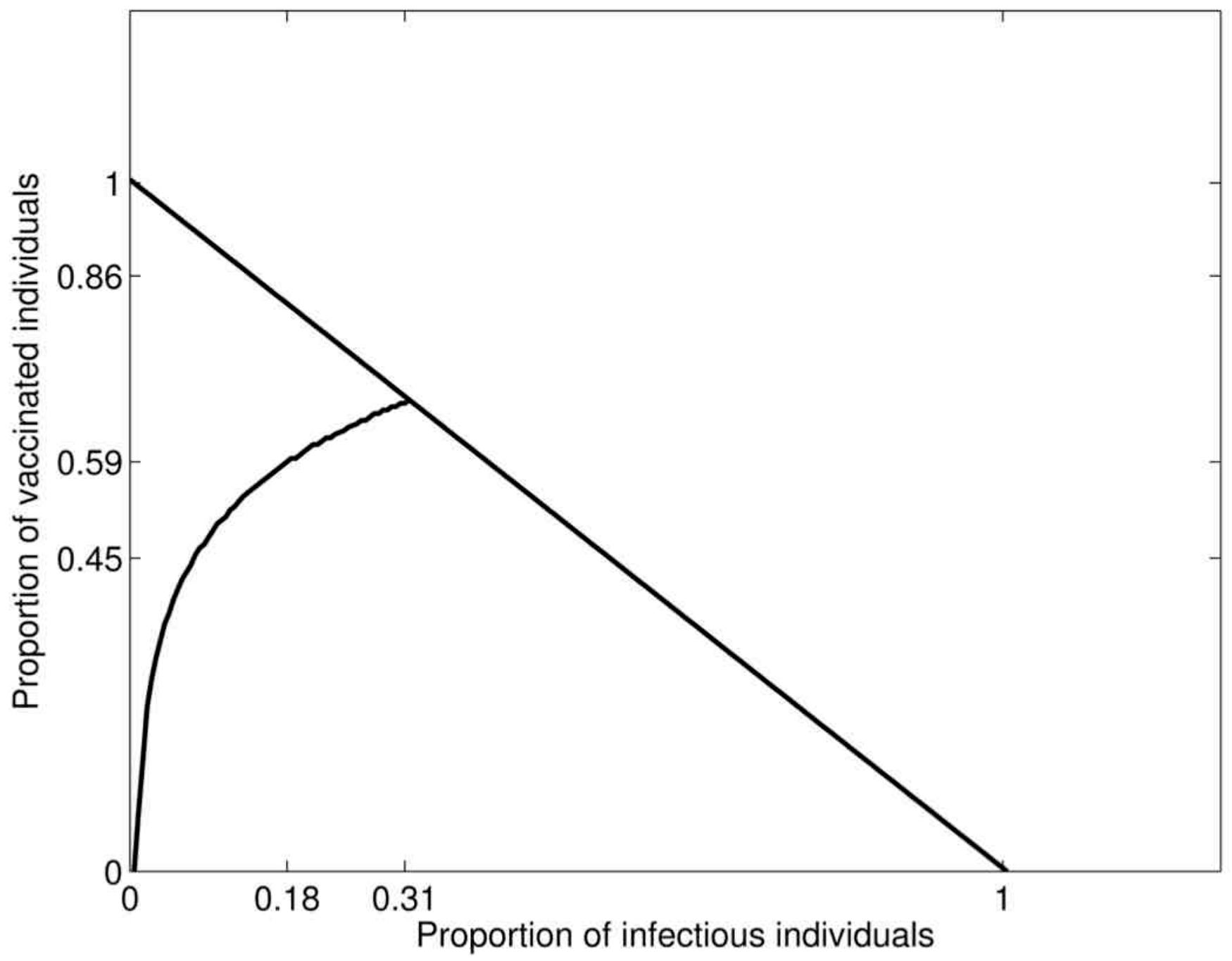}
\caption{The characteristic boundary separating the domains of attractions of $\bar x=(0,0.86)^\top$ (on the left) and $x^*=(0.31,0.45)^\top$ (on the right). The unstable equilibrium $\tilde x=(0.18,0.59)^\top$ is on the characteristic boundary. Note the slight abuse of notation due to the reduction from dimension three to two. The parameters are as in Figure~3.}
\end{figure}

\section{Computation of $\bar V$} \label{SecVBar}
In order to answer  the question of when the disease dies out, we have to compute the quantity $\bar V$ (cf.~equations~\eqref{EqExitTime} and~\eqref{EqVBar}). In other words, we have to solve a control problem. We first fix the time horizon $T>0$. We recall that we only have to consider absolutely continuous trajectories (as else $I_{T,x}(\phi)=\infty$) and hence we control the trajectory $\phi$ via its ``derivative''; for $t\in [0,T)$, $s \in [t,T]$ and an integrable function $\alpha$, we write
\[
\phi^\alpha(s):=\phi(s)=\phi(t)+\int_t^s \alpha(r) dr;
\]
for $x \in O$, we define the \emph{set of admissible controls} starting from $(t,x)$ by
\begin{equation} \label{EqVabrAdm}
\mathcal{A}^T(t,x):=\Big\{ \alpha:[t,T]\rightarrow \mathds{R}^d\text{ integrable } \Big| \; x+\int_t^s \alpha(r)dr \in \bar O \, \forall s\in [t,T] \text{ and } x+\int_t^T \alpha(r)dr \in A \setminus O \Big\}.
\end{equation}
The cost of an admissible control $\alpha$ starting from $(t,x)$ is given by
\[
J^T(t,x,\alpha)=\int_t^T L(\phi^\alpha(s),\alpha(s))ds,
\]
and thus the value function of the optimization problem becomes ($t\in [0,T)$)
\begin{equation}\label{EqValueFct}
v^T(t,x)=\inf_{\alpha \in \mathcal{A}^T(t,x)} J(t,x,\alpha)
\end{equation}
and
\[
v^T(T,x)=
\begin{cases}
0 & \text{if } x \not \in O \\
\infty & \text{else}.
\end{cases}
\]
The value function thus has a singularity at the terminal time $T$ which is due to the constraint $\phi(T)\not\in O$. For the models we have in mind (cf.~Examples~\ref{ExampleSIS} and~\ref{ExampleSIV}), it is possible to move on $\partial O \cap (A\setminus O)$ by following the solution of the ODE~\eqref{EqODE} (which is free of cost, cf.~\eqref{EqODE=0}).\footnote{For the SIS model, we have $\partial O \cap (A\setminus O)=\{0\}$. For the model with vaccination, we have
\[
\partial O \cap (A\setminus O)=\{ z\in A| \lim_{t \rightarrow \infty} Z^z(t)=\tilde x\}
\]
(recall that $\tilde x$ is the unstable endemic equilibrium).}
For this reason, we have $v^T(t,x) \leq v^{\tilde T}(t,x)$ for $\tilde T < T$ and therefore
\[
\lim_{T\rightarrow \infty} v^T(0,x^*)=\bar V.
\]
We want to remark here that any admissible control $\alpha$ with $\phi^\alpha(T) \not\in (A\setminus O) \cap \partial O$ (cf.~\eqref{EqODE=0} again) cannot be optimal. Therefore, it would be equivalent to require the weaker restriction $x+\int_0^s \alpha(r)dr \in A$ for $s\in[0,T]$ in equation~\eqref{EqVabrAdm}.

In the following, we describe a discrete dynamic programming algorithm for the numerical approximation of $\bar V$ (section~\ref{SubSecVbarDiscrete}) and apply this algorithm to the models from Example~\ref{ExampleSIS} (section~\ref{SubSecVBarSIS}) and Example~\ref{ExampleSIV} (section~\ref{SubSecVBarSIV}).

\subsection{Approximation by discrete-time dynamic programming} \label{SubSecVbarDiscrete}

In order to numerically compute the solution of the optimization problem~\eqref{EqValueFct}, we apply discrete-time dynamic programming. To this end, we discretize time and consider $n+1$ ($n \in \mathds{N}$) equidistant time points
\[
0=t_0<\dots< t_n=T, \quad \text{i.e.}, \quad t_m^n=t_m=m\Delta t \text{ for } \Delta t=\frac{T}{n}.
\]
For $j=0, \dots,n-1$ let $\tilde \alpha=(\tilde \alpha (t_m))_{m=j,\dots,n-1}$, where $\tilde \alpha(t_m) \in \mathds{R}^d$. We consider piecewise constant trajectories and define recursively for $x \in O$,
\[
\tilde \phi^{\tilde \alpha}(t_j)=\tilde \phi(t_j)=x, \quad \tilde \phi(t_{m+1})=\tilde \phi(t_m)+ \tilde \alpha(t_m) \Delta t \quad (m=j,\dots,n-1).
\]
We define the set of admissible controls by
\[
\tilde{\mathcal{A}}^{T,n}(t_j,x):=\Big\{ \tilde \alpha \; \Big|\; \tilde \phi^{\tilde \alpha}(t_m) \in A \, \forall m=j,\dots,n-1 \text{ and } \tilde \phi^{\tilde \alpha}(t_n) \in (A\setminus O) \cap \partial O \Big\}.
\]
The cost of an admissible control $\tilde \alpha \in \tilde{\mathcal{A}}^{T,n}(t_j,x)$ is given by
\[
\tilde J^{T,n}(t_j,x,\tilde \alpha):=\Delta t \sum_{m=j}^{n-1} L(\tilde \phi^{\tilde \alpha}(t_m),\tilde \alpha(t_m))
\]
and the value function is given by
\begin{equation} \label{EqOptDis}
\tilde v^{T,n}(t_j,x):=\inf_{\tilde \alpha \in \tilde{\mathcal{A}}^{T,n}(t_j,x)} \tilde J^{T,n}(t_j,x,\tilde \alpha).
\end{equation}
The constraint implies that
\begin{equation}\label{EqOptn-1}
\tilde v^{T,n}(t_{n-1},x)=\Delta t\inf_{\tilde \alpha \in \mathds{R}^d:x +\tilde \alpha \Delta t  \in (A\setminus O)\cap \partial O} L(x,\tilde \alpha).
\end{equation}
Furthermore, for $j<n-1$, the dynamic programming principle implies the following Bellman equation:
\begin{equation*}
\tilde v^{T,n}(t_j,x)=\inf_{\tilde\alpha \in \mathds{R}^d: x + \tilde \alpha\Delta t \in A} \Big\{
 L(x,\tilde \alpha) \Delta t + \tilde v^{T,n}(t_{j+1},x+\tilde \alpha\Delta t)
\Big\}.
\end{equation*}
Hence the solution of the optimization problem can be computed backward in time.

We need next to discretize space. Note also that in the SIS model, the functional $L$ is explicitly known, while this is not the case in the
vaccination model. For that reason we now distinguish between the two cases.

\subsection{The SIS model} \label{SubSecVBarSIS}

We discretize space via a grid of the set $A=[0,1]$  with fineness $\Delta x = 1/\bar n$, for some new integer $\bar n$. Hence the grid points are given by
\[
x^i :=i \Delta x \quad \text{for } i=0,\dots,\bar n.
\]

Given a function $f$ which is only defined on the grid points, we define values of the function for arbitrary $x \in A$ by its linear interpolation $f_{\inter}$. In other words, we obtain for $x\in A$,
\[
f_{\inter}(x):= \bar n \sum_{i=0}^{\bar n} \mathds{1}_{\{x \in [x^{i}, x^{i+1})\}}
 \Big\{ f(x^{i}) (x^{i+1}-x) +f(x^{i+1}) (x-x^{i}) \Big\}.
\]

We obtain the following (backward in time) algorithm for the solution of the optimization problem~\eqref{EqOptDis}.
\begin{enumerate}
\item[1.]
For $x$ on the grid,
\[
\bar v^{T,n,\bar n}(t_{n-1},x) = \Delta t\inf_{\tilde \alpha \in \mathds{R}^d:x + \tilde \alpha\Delta t \in (A\setminus O)\cap \partial O} L(x,\tilde \alpha).
\]
We denote the solution of this minimization problem by $\tilde \alpha^*(t_{n-1})=\tilde \alpha^*(t_{n-1},x)$.
\item[2.]
For $m < n-1$, $x$ on the grid,
\[
\bar v^{T,n,\bar n}(t_{m},x) = \inf_{\tilde\alpha \in \mathds{R}^d: x + \tilde \alpha\Delta t  \in A} \Big\{ L(x,\tilde \alpha) \Delta t + \bar v_{\inter}^{T,n,\bar n}(t_{m+1}, x+\tilde \alpha\Delta t)
\Big\}.
\]
We denote the solution of this minimization problem by $\tilde \alpha^*(t_m)=\tilde \alpha^*(t_m,x)$.
\end{enumerate}
Note that this algorithm involves the solution of $n\times \bar n$ minimization problems. We recover an approximation $\tilde \phi^*=(\phi^*(t_m))_{m=0,\dots,n}$ of the optimal trajectory  starting from the stable equilibrium $x^*$, by the following procedure (including another linear interpolation).
\begin{enumerate}
\item[3.]
For $m=0,\dots,n-1$,
\[
\tilde \phi^*(t_0)=x^*, \quad \tilde \phi^*(t_{m+1})=\phi^*(t_m)+ \tilde \alpha^*_{\inter}(t_m,\phi^*(t_m)) \Delta t.
\]
\end{enumerate}

It is well-known that under appropriate assumptions, the solution of this procedure converges to the solution of the optimization problem~\eqref{EqValueFct} (see, e.g., \cite{Falcone1994,Falcone2013}); the optimal coupling of time and space discretization is $\Delta t=\Delta x$. We want to remark here that standard assumptions do not apply to our model. For instance, the function $L$ is not Lipschitz continuous, even $L(x,y)\rightarrow \infty$ for $x \rightarrow \partial A$ is possible.\footnote{This is e.g. the case for $y<0$ and $x\rightarrow 0$. Note however that the action required to approach the boundary, e.g., at constant speed one is finite and we have $\int_0^t L(t-s,-1) ds \rightarrow 0$ as $t \rightarrow 0$.} Furthermore, the terminal constraint results in a singularity of the value function at terminal time $T$. We are not aware of any analytical results including our set--up. We nevertheless have reasons to believe that
\begin{equation*}
\bar V \approx \bar v^{T,n,\bar n}_{\inter}(t_0,x^*) \quad \text{for large enough } T.
\end{equation*}
It is not in the scope of this article to provide a thorough analytical proof of the convergence of the algorithm.

We readily observe that for the SIS model, the Legendre-Fenchel transform is given explicitly by
\[
L(x,y)=y \log \tilde \theta - \beta x (1-x) (\tilde \theta - 1) - \gamma x \Big(\frac{1}{\tilde \theta}-1 \Big),
\]
where
\[
\tilde \theta:=\frac{y+\sqrt{y^2+4 \beta\gamma x^2 (1-x)}}{2\beta x(1-x)}.
\]
Furthermore, the boundary consists of the single point $\{0\}$; hence equation~\eqref{EqOptn-1} becomes
\[
\tilde v^{T,n}(t_{n-1},x)=\Delta t L(x,-x/\Delta t).
\]
We consider the following parameters: $\beta=1.5$, $\gamma=1$; in particular, we have $R_0=1.5>1$ and the endemic equilibrium $x^*=1/3$ is stable. We apply the dynamic programming algorithm described in section~\ref{SubSecVbarDiscrete} with $\Delta t=\Delta x =1/100$ and $T=5,10,20,40,60$ (i.e., $\bar n=100$ and $n=500, 1000, 2000, 4000, 6000$). We obtain
\begin{align*}
v^{5,n,\bar n}_{\inter}(t_0,x^*)& =0.0953, \quad v^{10,n,\bar n}_{\inter}(t_0,x^*)=0.0757, \quad v^{20,n,\bar n}_{\inter}(t_0,x^*)=0.0705, \\
v^{40,n,\bar n}_{\inter}(t_0,x^*)&=0.0702, \quad v^{60,n,\bar n}_{\inter}(t_0,x^*) =0.0702.
\end{align*}
We hence deduce that $\bar V\approx 0.0702$ and thus $\mathbb{E}[\tau^{N,x}]\approx \exp(0.0702 N)$ for large $N$. We illustrate the solution of the optimization problem in Figure~5 for $T=5, 20, 60$. We observe that it is the cheapest to move slowly near the stable equilibria $x^*$ and the unstable equilibrium $\bar x=0$.

\begin{figure} 
\centering \includegraphics[height=0.225\hsize]{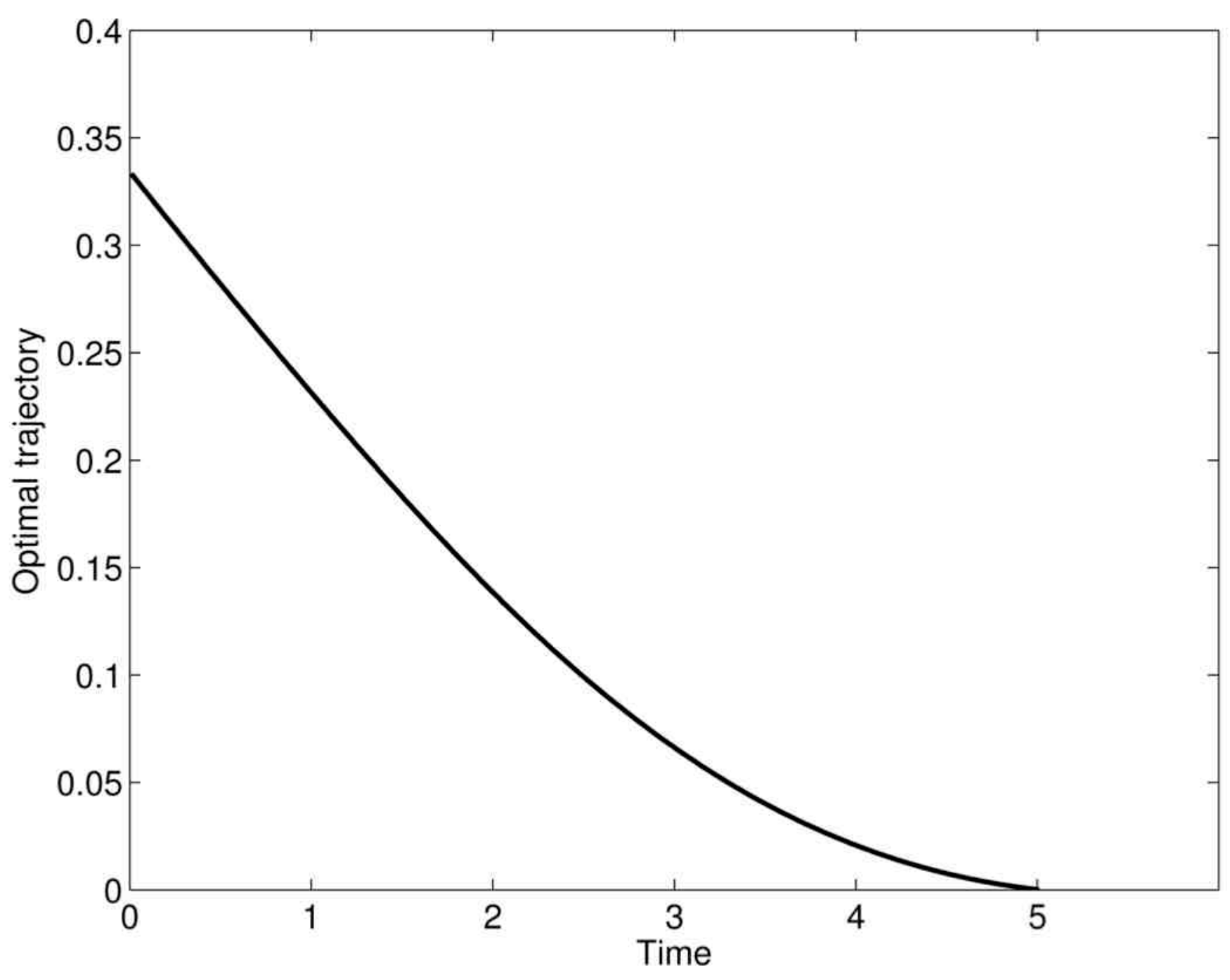}
\quad
\includegraphics[height=0.225\hsize]{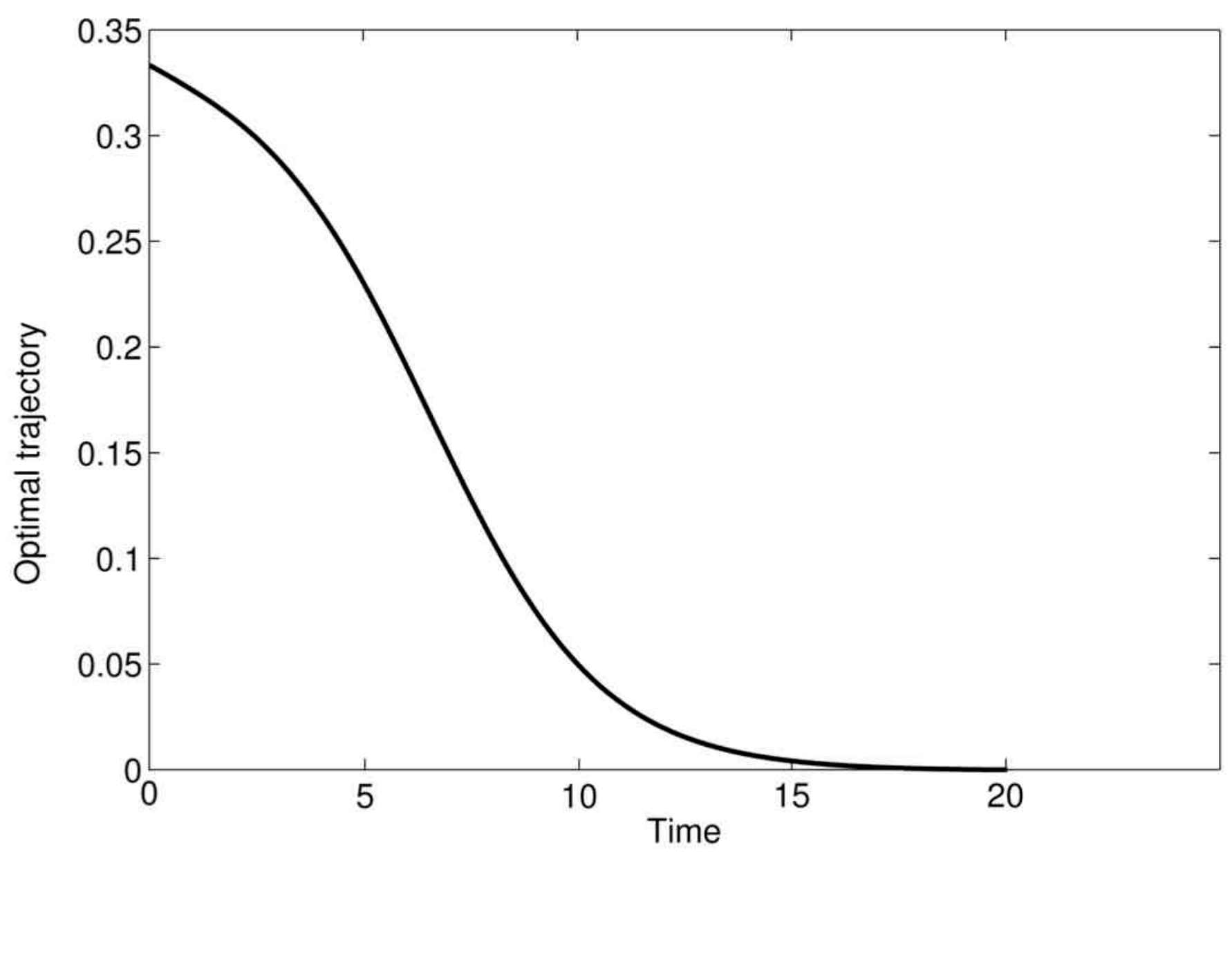}
\quad
\includegraphics[height=0.225\hsize]{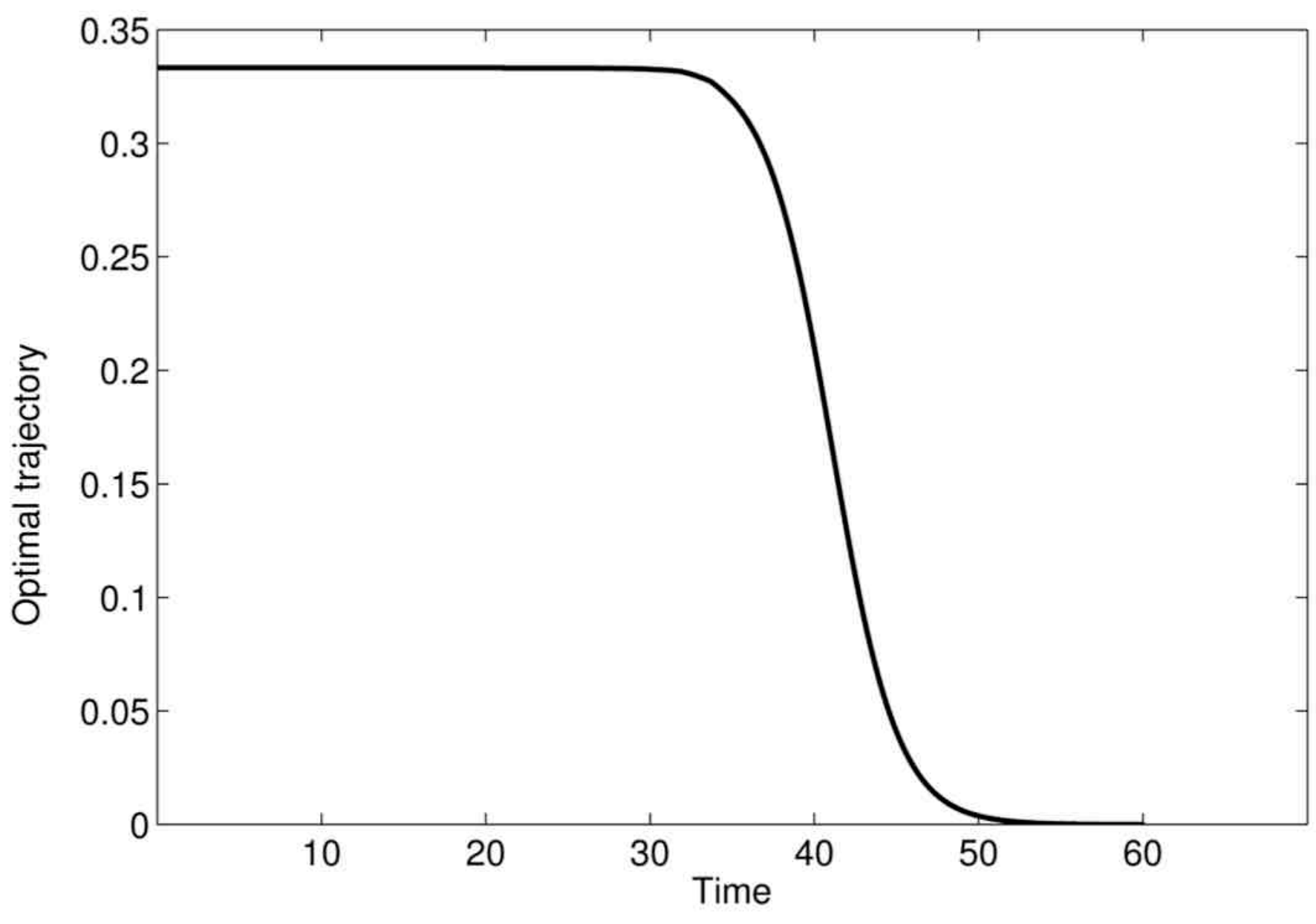}
\caption{Optimal trajectory $\tilde \phi^*$ for $T=5$ (left picture), $T=20$ (middle picture) and $T=60$ (right picture); $\Delta t=\Delta x=\frac{1}{100}$, $\beta=1.5$ and $\gamma=1$.}
\end{figure}

\subsection{The model with vaccination} \label{SubSecVBarSIV}

Unlike in section~\ref{SubSecVBarSIS}, the place of exit from $O$ is a priori not clear for the model with vaccination (cf.~Example~\ref{ExampleSIV}). It seems
plausible that it is cheapest to exit $O$ at the unstable endemic equilibrium $\tilde x$. Indeed, if $T$ is chosen large enough, an optimal (or near optimal) trajectory $\phi$ can move on the boundary towards $\tilde x$ after the exit time without generating any extra cost (cf.~\eqref{EqODE=0}); in this case we have $\lim_{t \rightarrow \infty} \phi(t)=\tilde x$. Hence, we can assume that $\tilde x$ is indeed the place of exit for the numerical computation of $\bar V$. Therefore, equation~\eqref{EqOptn-1} becomes
\begin{equation*}
\tilde v^{T,n}(t_{n-1},x)=\Delta t L(x,(\tilde x- x)/\Delta t).
\end{equation*}

However, the Legendre-Fenchel transform $L$ cannot be computed explicitly. Here we use the following definition of $L$, which is equivalent to \eqref{EqLegendre} (see \cite{Kratz2014}):
\begin{equation*}
L(x,y):=\begin{cases}
\inf_{\mu \in B_{x,y}} \tilde \ell (\mu,x)& \text{if } B_{x,y} \not=\emptyset \\
\infty& \text{if not},
\end{cases}
\end{equation*}
where for $\mu \in B_{x,y}$,
\begin{align*}
\tilde \ell(\mu,x)&:= \sum_{j} \big\{ \beta_j(x) - \mu_j + \mu_j \log\big(\tfrac{\mu_j}{\beta_j(x)}\big)\big\};
\end{align*}
and
\begin{align*}
B_{x,y}&:= \Big\{\mu \in \mathds{R}^k_+ | \mu_j >0 \text{ only if } \beta_j(x)>0 \text{ and } y=\sum_j \mu_j h_j\Big\}.
\end{align*}
Since we seek to minimize
$$\overline V:=\inf_{T, \phi\text{ s.t. }\phi(0)=x^*, \phi(T)=\tilde x}\int_0^T L(\phi(t),\phi'(t))dt,$$
our discretized dynamic programming equation reads: for each $x\in\mathcal{G}$ (the space grid),
\begin{align*}
v(t_k,x)&=\inf_{x'\in\mathcal{G}}\Big\{L\left(x,\frac{x'-x}{\Delta t}\right)\Delta t+v(t_{k+1},x')\Big\} \\
&=\inf_{x'\in\mathcal{G}}\Big\{\tilde \ell(\hat\mu(x,x'),x)\Delta t+v(t_{k+1},x')\Big\},
\end{align*}
where $\hat\mu(x,x')$ is the argument which realizes the following minimum
$$\inf_{\mu,\ x+\Delta\sum_j\mu_jh_j=x'}\tilde\ell( \mu,x).$$
Note that here we avoid to consider interpolation between points on the grid. We have implemented this algorithm
on a grid whose one of the axis is the line passing through  $x^*$ and $\tilde x$, with the following values of the parameters :
$\beta=3.6$, $\theta=0.02$, $\mu=0.03$, $\eta=0.3$, $\sigma=0.1$ and $\gamma=1$.

The approximation of $\bar V$ which we have found is $0.3891$, with $\Delta t=0.05$, and the $\Delta x_1=\Delta x_2=0.005$.
This means that $\mathbb{E}[\tau^{N,x^*}]\approx \exp(0.3891 N)$ for large $N$.
Note that if $\Delta t$ is not small enough compared to $\Delta x$, since we restrict the above minimization to $x'\in\mathcal{G}$,
we cannot allow small speed, and we get larger values of $\bar V$. We expect that more accurate results could be obtained,
by doing the minimization over $x'$ on a finer grid, which necessitates to interpolate $v(t_{k+1},\cdot)$ between points on the grid.
We have not been able to implement this refinement, by lack of time. This will be done in the future. It should improve the numerical
estimate of $\bar V$. The above quantity should be considered as an upper bound of the true value.

\section{Simulation of $Z^N$} \label{SecSim}
In this section, we describe known exact and approximate simulation methods for the simulation of the solution $Z^N$ of equation~\eqref{EqPoisson}. To this end, we fix $N$; with the notations $\nu_{j}=\frac{h_{j}}{N}$ and $a_{j}(z)=N\beta_{j}(z)$, equation~\eqref{EqPoisson} can be rewritten as
\begin{equation}\label{EqPoissonMod}
  Z^{N}(t) =x+\sum_{j=1}^{k}\nu_{j}P_{j}\Big(\int_{0}^{t}a_{j}(Z^{N}(s))ds\Big).
\end{equation}

Note that our ultimate goal would be to simulate the value of $\tau^N$ which, from the above Large Deviations results, is huge, as soon as $N$ is large.
This may be seen as impossible. However, this is similar to the problem of rare event simulation,  for which there exist specific techniques, which we intend to
implement in the future. See also the discussion below in section \ref{SubSubSecSimRare}. In the framework of this paper, we present some preliminary numerical work, which aims only at simulating the stochastic process
$Z^N$ over a given finite time horizon.

In the following, we first introduce the so-called stochastic simulation algorithm (SSA) proposed in~\cite{Gillespie1976, Gillespie1977} (section~\ref{SubSecSimSSA}). The SSA allows for the exact simulation of the process $Z^N$. Unfortunately, this algorithm is rather slow in some circumstances, typically when  there are too many jumps of the Poisson processes on a fixed time interval. Therefore, we discuss faster approximate simulation methods in the subsequent sections: in section~\ref{SubSecSimTauGeneral}, we introduce the explicit $\tau$-leaping algorithm (see~\cite{Gillespie2001}) in the form of~\cite{Cao2003}; in section~\ref{SubSecSimTauMod}, we discuss modifications of the $\tau$-leaping algorithm: implicit and mid-point $\tau$-leaping (section~\ref{SubSubSecSimTauImp}) and the algorithm from~\cite{Cao2005} which ensures that the simulated process remains in its domain, e.g., by preventing the components of the process from becoming negative (section~\ref{SubSubSecSimTauMod}).
We apply the SSA and the modified $\tau$-leaping algorithm from section~\ref{SubSubSecSimTauImp} to the specific models from Examples~\ref{ExampleSIS} and~\ref{ExampleSIV} in sections~\ref{SubSecSimSIS} and~\ref{SubSecSimSIV}, respectively. We conclude by addressing open research questions concerning the simulation of the process $Z^N$ in section~\ref{SubSecSimOpen}.

\subsection{The stochastic simulation algorithm} \label{SubSecSimSSA}

A mathematically exact procedure for simulating the evolution of the process $Z^N$ given by equation \eqref{EqPoissonMod} is the stochastic simulation algorithm (SSA) first proposed in \cite{Gillespie1976, Gillespie1977}. The SSA changes the value of the solution at each jump time of one of the Poisson processes. The simplest implementation of the SSA is the so-called \emph{direct method}:

Whenever $t<T$:
\begin{enumerate}
\item[1.]
While in state $z$ at time $t$, evaluate all functions $a_j(z)$ and their sum $a_{0}(z):=\sum_{j=1}^{k}a_{j}(z)$.
\item[2.]
Generate the time increment $\tau$ as the realization of an exponential random variable with parameter $a_{0}(z)$.
\item[3.]
Choose an index $j\in \{1,\dots,k\}$ according to the probability distribution $\frac{a_{j}(z)}{a_{0}(z)}$ ($j=1,\dots,k$).
\item[4.]
Update
\begin{enumerate}
\item[i.]
$t\leftarrow t+\tau$,
\item[ii.]
$z\leftarrow z+\nu_{j}$.
\end{enumerate}
 \item[5.]
 Record $(t,z)$. Return to 1.~if $t+\tau<T$; else stop.
\end{enumerate}

Carrying out step 3.~is mathematically straightforward: we draw one realization $r$ of a uniform random variable on $[0,1]$; then, $j$ is the smallest positive integer
for which
\begin{equation*}
  \sum_{i=1}^{j}a_{i}(z)>ra_{0}(z).
\end{equation*}
The SSA is mathematically exact. However, the task of explicitly simulating each jump makes the SSA too slow for practical implementation in some circumstances. The reason is that for large $N$, there are too many jumps on any fixed time interval.

\subsection{The explicit Poisson $\tau$-leaping algorithm} \label{SubSecSimTauGeneral}

A faster but approximate stochastic simulation procedure is the explicit Poisson $\tau$-leaping algorithm (see~\cite{Gillespie2001}). The basic idea of this procedure is to advance the system by a \emph{preselected} time increment $\tau$ (in
contrast to the \emph{randomly generated} time increment $\tau$ in the SSA). On the one hand, $\tau$ should be chosen large enough in order to ensure that ``many'' jumps occur during this length of time (as else, the savings in computation time are not significant); on the other hand, it should be small enough such that the values of the functions $a_j$ are unlikely to change ``significantly'' (as else, the procedure is not accurate enough). The latter restriction is called the \emph{leap condition}. A strategy due to~\cite{Cao2003} for satisfying this condition is to require that the changes of the functions $a_j$ during a leap $\tau$ are ``likely'' to be bounded by $\epsilon a_{0}(z)$; here, $\epsilon$ ($0<\epsilon\ll1$) is the error control parameter. In other words,
\begin{equation*} 
|a_{j}(Z^{N}(t+\tau))-a_{j}(Z^{N}(t))|\leq\epsilon a_{0}(Z^{N}(t)) \quad \text{for }  j=1\dots,k \text{ with high probability.}
\end{equation*}

We estimate the largest value of $\tau$ that meets this particular requirement as follows. First, we compute
\begin{equation}\label{EqTauSel1}
  f_{jj'}(z):=\sum_{i=1}^{d}\frac{\partial a_{j}(z)}{\partial z_{i}}\nu_{ij'}\quad \text{for }j,j'=1,\dots,k;
\end{equation}
here $\nu_{ij}$ denotes the $i$-th component of the vector $\nu_j$. Subsequently, we calculate
\begin{equation}\label{EqTauSel2}
  \mu_{j}(z):=\sum_{j'=1}^{k}f_{jj'}(z)a_{j}(z) \quad \text{and} \quad
  \sigma_{j}^{2}(z):=\sum_{j'=1}^{k}f_{jj'}^{2}(x)a_{j}(z) \quad \text{for } j=1, \dots,k
\end{equation}
and set
\begin{equation}\label{EqTauSel3}
  \tau=\tau(z):=\min_{j=,1\dots,k}\left\{\frac{\epsilon a_{0}(z)}{|\mu_{j}(z)|},\frac{\epsilon^{2} a_{0}^{2}(z)}{\sigma_{j}^{2}(z)}\right\}.
\end{equation}
Note that  for a standard Poisson process $P$,
\begin{align*}
 P_{j}\Big(\int_{0}^{t+\tau}a_{j}(Z^{N}(s))ds\Big)-P_{j}\Big(\int_{0}^{t}a_{j}(Z^{N}(s))ds\Big) &\stackrel{d}{=}P\left(\int_{t}^{t+\tau}a_{j}(Z^{N}(s))ds\right) \\
   &\approx P(a_{j}(Z^{N}(t))\tau) = P(a_{j}(z)\tau)=:p_{j}
\end{align*}
and thus
\begin{equation*}
  Z^{N}(t+\tau) \stackrel{d}{\approx}  Z^{N}(t)+\sum_{j=1}^{k}p_{j}\nu_{j}
                =  z+\sum_{j=1}^{k}p_{j}\nu_{j}.
\end{equation*}
This procedure for the selection of $\tau$ yields the following $\tau$-leaping algorithm.

Whenever $t<T$:
\begin{enumerate}
\item[1.]
While in state $z$ at time $t$, evaluate the functions $a_j$ and their sum $a_{0}(z)=\sum_{j=1}^{k}a_{j}(z)$.
\item[2.]
Compute $\tau$ via the formulas~\eqref{EqTauSel1}-\eqref{EqTauSel3}.
\item[3.]
Fix a small integer $n$ (typically, e.g., $n=10$) and a moderate integer $\bar n$ (typically, e.g, $\bar n=100$).
\begin{enumerate}
\item[i.]
If $\tau < \frac{n}{a_{0}(z)}$, reject $\tau$ and execute $\bar n$ iterations of the SSA. Then, return to step 1.
\item[ii.]
If $\tau\geq\frac{n}{a_{0}(z)}$, proceed to step 4.
\end{enumerate}
\item[4.]
For $j=1,\dots,k$, generate $p_{j}$ as a   Poisson random variable with mean $a_{j}(z)\tau$.
\item[5.] Update
\begin{enumerate}
\item[i.]
$t\leftarrow t+\tau$,
\item[ii.]
$z\leftarrow z+\sum_{j=1}^{k}p_{j}\nu_{j}$.
\end{enumerate}
 \item[6.]
 Record $(t,z)$. Return to 1.~if $t+\tau<T$; else stop.
\end{enumerate}

Let us shortly comment on this algorithm. $p_{j}$ represents an approximation of the number of jumps of type $j$ during the time interval $[t,t+\tau[$. We remark again that the approximation by $p_{j}$ in step 4.~is justified as long as $a_{j}(z)$ remains more or less constant over the subsequent time interval of length $\tau$. The reason for step 3.~i.~is the following. $\frac{1}{a_{0}(z)}$ is the mean time step to the next jump when applying the SSA. If $\tau$ is not significantly larger than $\frac{1}{a_{0}(z)}$, the computation of $\tau$ (via equations~\eqref{EqTauSel1}-\eqref{EqTauSel3}) is  no more efficient then applying the SSA directly (which is also of course more accurate); as this property is likely to persist for a while, it is more appropriate to apply the SSA for some time before returning to the $\tau$-selection procedure again.

The explicit $\tau$-leaping algorithm above has been shown to provide accurate simulation results; furthermore, it is substantially faster than the SSA for many ``not-too-stiff'' systems\footnote{I.e., systems for which the difference between the smallest and the largest jump rates is not too large (see, e.g., \cite{Cao2005}).\label{FootnoteStiff}} (see~\cite{Cao2005}; see also~\cite{Anderson2011} for theoretical evidence for this). We illustrate the gain of computation time by numerical experiments in sections~\ref{SubSecSimSIS} and~\ref{SubSecSimSIV} below for a modification of this algorithm (cf.~section~\ref{SubSubSecSimTauMod}).

\subsection{Modifications of the $\tau$-leaping algorithm} \label{SubSecSimTauMod}
There exist several variants of the above described explicit $\tau$--leaping algorithm. We now describe some of these.

\subsubsection{Implicit and mid--point $\tau$-leaping algorithms} \label{SubSubSecSimTauImp}
One may want to freeze the rates not at the value of the solution at the beginning of the discretization interval, but rather
at an approximate value of the position at the end of that interval, or else at the mid-point.

The ``implicit'' $\tau$-leaping algorithm is obtained by replacing $a_j(Z^N(t))\tau$ with
\[
a_j\big(Z^N(t)+\tau b(Z^N(t))\big) \tau,
\]
and the mid--point $\tau$-leaping algorithm consists in replacing the same quantities by
\[
a_j\big(Z^N(t)+\tau b(Z^N(t))/2\big) \tau.
\]
In~\cite{Anderson2011} it is proved that the error in the mid-point $\tau$-leaping algorithm is less than the one produced by the
explicit $\tau$--leaping algorithm.

\subsubsection{A modified $\tau$-leaping algorithm} \label{SubSubSecSimTauMod}
As most numerical schemes, the $\tau$-leaping algorithm has the disadvantage
that the approximate solution might very well jump over the boundaries which the solution  of the stochastic equation~\eqref{EqPoissonMod} can never cross. More precisely, the solution of~\eqref{EqPoisson} satisfies certain constraints. For example, in the case of constant total population, we have that $Z^N_i(t)\ge0$, $1\le i\le d$, and $Z^N_0(t):=1-\sum_{i=1}^dZ^N_i(t)\ge0$.
We would like that the approximate solution satisfies these constraints as well. However, this might not be the case with
a standard $\tau$-leaping algorithm. In particular, if for some $0\le i\le d$, $Z^N_i(t)$ is very small, then some of the Poisson driving processes may push it over to a negative value during the next time interval of length $\tau$.
Because of that, at the beginning $t$ of each time step, we define
\[
L_j:=\min_{0\le i\le d,\; \nu_{ij} a_j(Z^N(t))<0} Z^N_i(t),\quad
\text{where }  \nu_{0j}:=-\sum_{i=1}^d\nu_{ij},
\]
and choose a critical value $n_c$. The critical directions will be chosen by comparing $L_j$ to $n_c$.

The modified Poisson $\tau$-leaping procedure requires that  a critical Poisson process cannot jump more than once in a time interval of length $\tau$. That makes it impossible for a critical Poisson process to drive any component of $Z^{N}$ to a negative value. The theoretical justification for each step in this modified $\tau$-leaping procedure is given in~\cite{Cao2005}.
\par Whenever $t<T$:
\begin{enumerate}
\item[1.]
While in state $z$ at time $t$, evaluate the functions $a_j(z)$ and their sum $a_{0}(z)=\sum_{j=1}^{k}a_{j}(z)$.
\item[2.]
Identify the currently critical Poisson processes as those $P_{j}$ for which
\[
a_{j}(z)>0 \quad \text{and} \quad L_{j}<n_{c}.
\]
\item[3.]
Modify the $\tau$-selection procedure from equations~\eqref{EqTauSel1}-\eqref{EqTauSel3} by only taking into account those jump directions $j'$ corresponding to the \emph{non-critical} Poisson processes. Denote the resulting time step by $\tau'$. If there is no non-critical Poisson process, set $\tau'=\infty$.
\item [4.]
Fix a small integer $n$ (typically, e.g.,$n=10$) and moderate integer $\bar{n}$ (typically, e.g, $\bar{n}=100$).
\begin{enumerate}
\item[i.]
If $\tau'<\frac{n}{a_{0}(z)}$, reject $\tau'$ and execute $\bar{n}$ iterations of the SSA. Then, return to step 1.
\item[ii.]
If $\tau'\geq\frac{n}{a_{0}(z)}$ proceed to step 5.
\end{enumerate}
\item[5.]
Compute the sum $a_{0}^{c}(z)$ of the functions $a_{j}(z)$ of the critical Poisson processes. Generate $\tau''$ as a realization of the exponential random variable with parameter $a_{0}^{c}(z)$.
\item[6.]
\begin{enumerate}
\item[i.]
If $\tau'<\tau''$, set $\tau=\tau'$. For all the critical Poisson processes $P_{j}$, set $p_{j}=0$. For all the non-critical Poisson processes $P_{j}$, generate $p_{j}$ as a realization of the Poisson random variable with mean $a_{j}(z)\tau$.
\item[ii.] If $\tau''\leq\tau'$, set $\tau=\tau''$. Generate an integer $j_{c}$ according to the probability distribution $a_{j}(z) /a_{0}^{c}(z)$ (i.e., $j$ runs over the index values of the critical Poisson processes). Set $p_{j_{c}}=1$, and for all the other critical Poisson process, set $p_{j}=0$. For all the non-critical Poisson processes $P_{j}$, generate $p_{j}$ as a realization of the Poisson random variable with mean $a_{j}(z)\tau$.
\end{enumerate}
\item[7.]
Update
\begin{enumerate}
\item[i.] $t\leftarrow t+\tau$,
\item[ii.] $z\leftarrow z+\sum_{j=1}^{d}p_{j}\nu_{j}$.
\end{enumerate}
\item[8.]
If $z$ is not in its domain, replace $\tau'\leftarrow \tau'/2$, and return to step 6.
\item[9.]
Record (t,z). Return to 1 if $t+\tau<T$; else stop.
\end{enumerate}

\begin{rmrk}
It may seem strange to use a uniform critical value $n_c$, independent of the value of the rate of each Poisson process.
It would be natural to compare each $L_j$ with the corresponding $a_j(Z^N(t))\tau$, which however would require to compute a first value of $\tau$. We intend to search for improvements in that direction in the future.
\end{rmrk}

\begin{rmrk}
Although the modified $\tau$--leaping algorithm is designed so as to reduce the chance of a negative value of
any of the $Z^N_i(t)$, it cannot completely rule them out. This is the reason for step 8.~of the algorithm. We suspect that this step
might introduce a bias in the simulation, in that it tends to delete simulations with exceptionally large values of some of the Poisson process increments, and also to reduce the chance of approaching or hitting any of the boundaries, which is specially annoying for our purposes.
\end{rmrk}

We intend to check two other possible approaches. The first one consists in projecting the computed point, if it is outside the domain. The second consists in replacing $\tau$ by $\tau/2$, but without resimulating a new independent value of the Poisson process increment, but rather by simulating the new increments according to the conditional law of $P_j(t+\tau')-P_j(t)$,
given the previously simulated value of $P_j(t+\tau)-P_j(t)$, which, if $P_j(t+\tau)-P_j(t)=k$, follows a binomial law with parameters
$k$ and $\tau'/\tau$.

\subsection{The SIS model}\label{SubSecSimSIS}

In this section, we apply the simulation algorithms from sections~\ref{SubSecSimSSA} and~\ref{SubSubSecSimTauMod} to the SIS model (Example~\ref{ExampleSIS}). We consider an endemic situation ($\beta=1.5$, $\gamma=1$, hence $R_0=1.5>1$ and $x^*=1/3$) and compare the performance of the SSA to the modified $\tau$-leaping algorithm. If the population size is relatively small, $N=2000$, we observe that the computation times are similar for both algorithms (even smaller for the SSA). A closer look at the realized trajectory of the $\tau$-leaping procedure confirms that almost all jumps are of size $1/2000$. This suggests that in most cases, the selected time step $\tau'$ is rejected and the SSA is applied (cf.~steps 3.~and 4.~i.~of the algorithm). Hence, there is no significant gain in computation time. If we let $N=20000$, the number of jumps (and hence the computation time) for the SSA increases approximately by a factor ten. The $\tau$-leaping algorithm is significantly faster which suggests that the time step $\tau'$ is now rarely refused. While the SSA becomes inefficient for very large time horizons (e.g., $T=600$) and population sizes (e.g., $N=200000$), the process can still be simulated efficiently by the $\tau$-leaping algorithm.

We illustrate realized trajectories of the SSA and the modified $\tau$-leaping algorithm in Figures~6 and~7, respectively (left pictures: $N=2000$, right pictures: $N=20000$). We observe that they approach the solution of the ODE as $N$ increases for both algorithms as predicted by the LLN (cf.~\eqref{EqLLN}; see~\cite{Kurtz1978}).

\begin{figure} 
\centering \includegraphics[width=0.45\hsize,height=0.35\hsize]{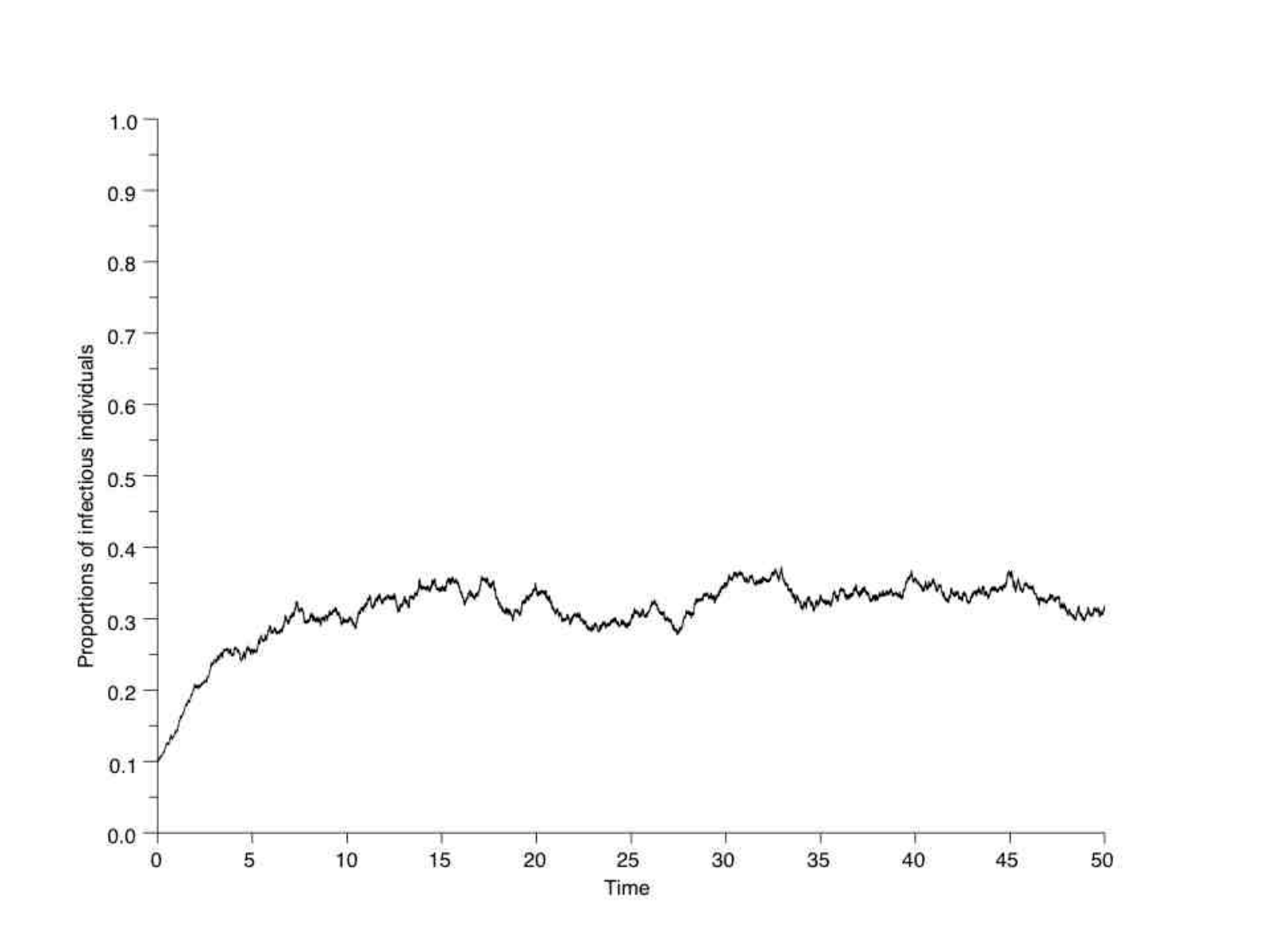}
\qquad
\includegraphics[width=0.45\hsize,height=0.35\hsize]{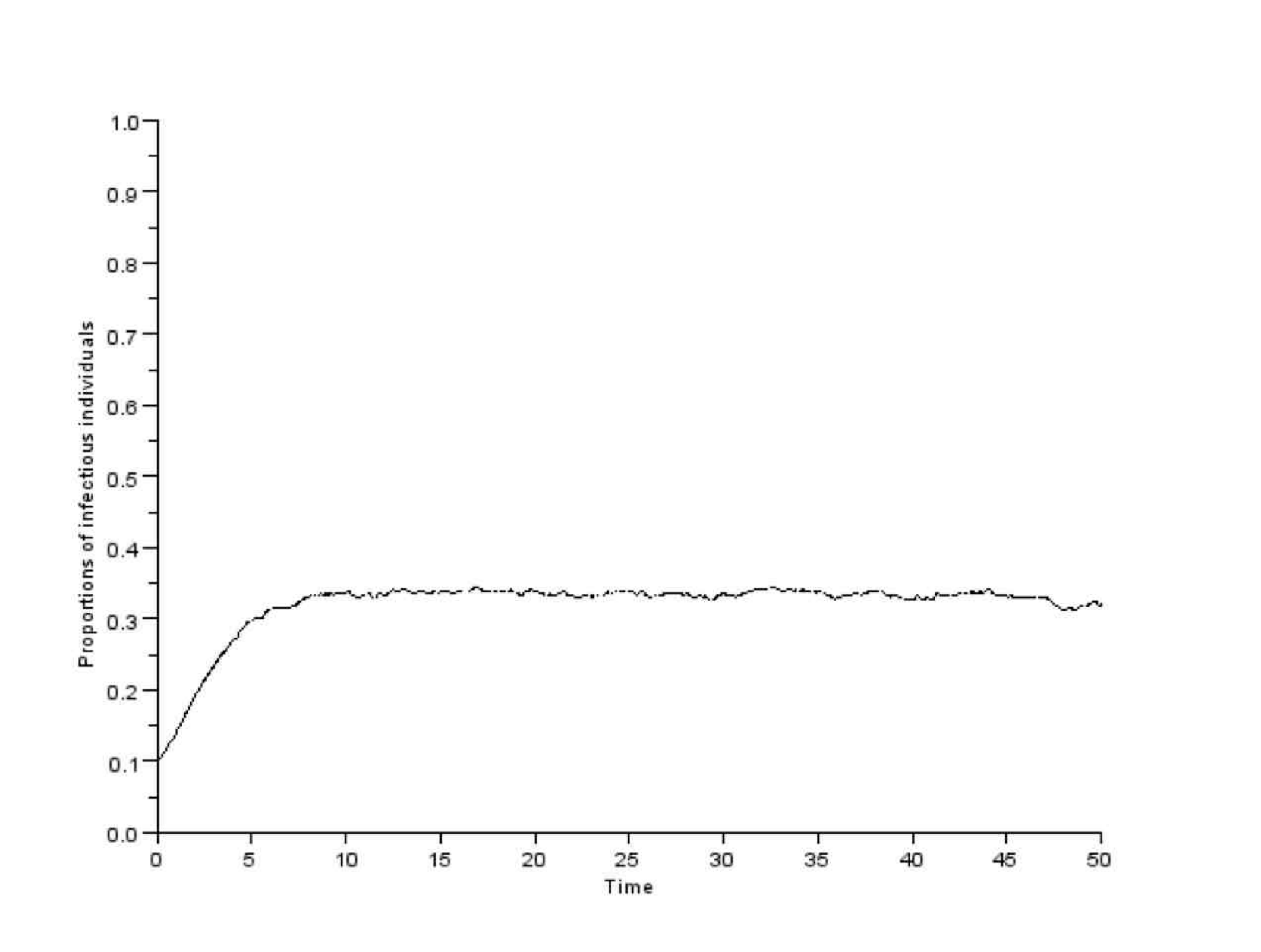}
\caption{Realized trajectories of the SSA for $N=2000$ (left picture) and $N=20000$ (right picture). $Z(0)=0.1$, $T=50$, $\beta=1.5$, $\gamma=1$.}
\end{figure}

\begin{figure} 
\centering \includegraphics[width=0.45\hsize,height=0.35\hsize]{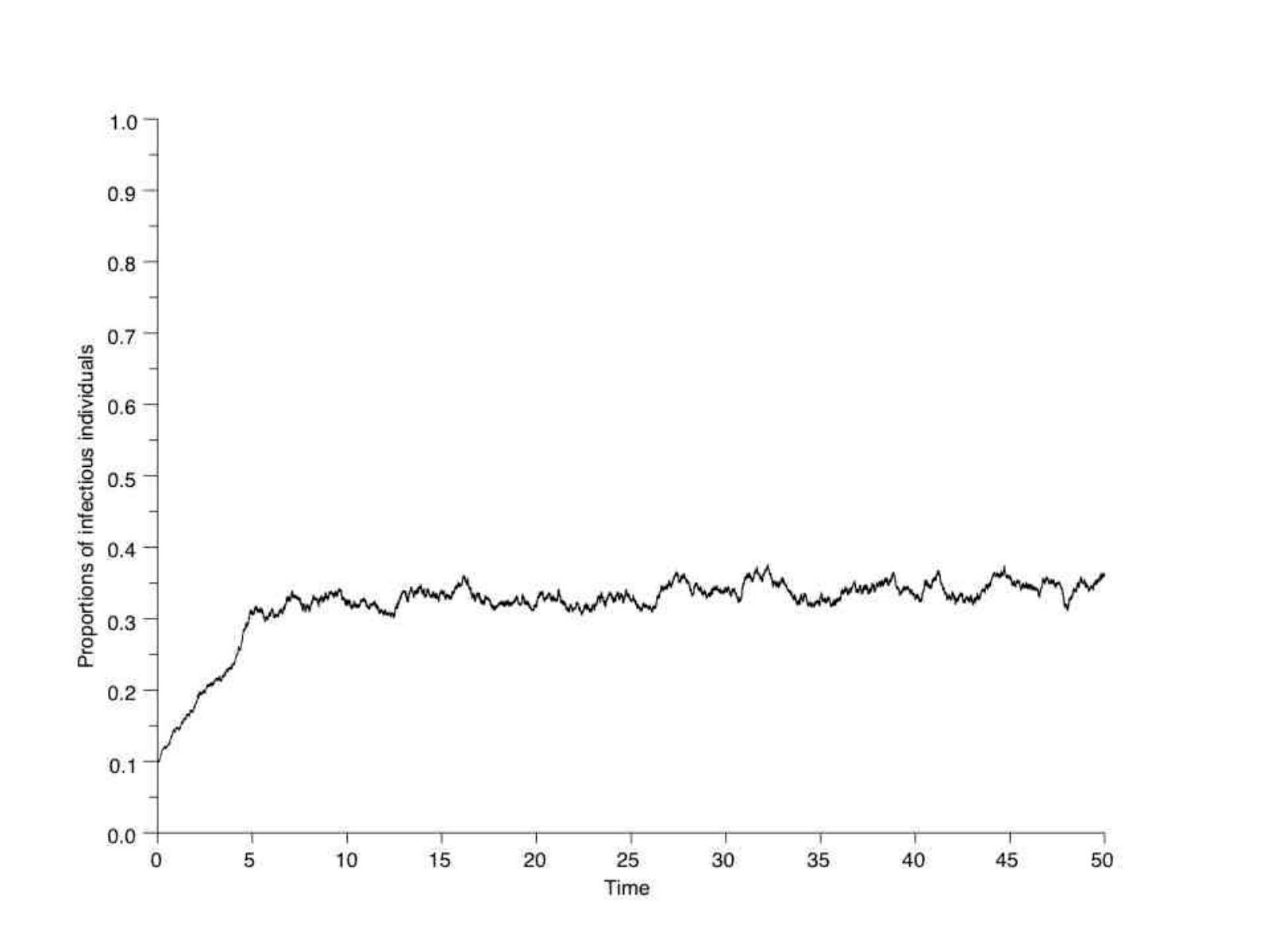}
\qquad
\centering \includegraphics[width=0.45\hsize,height=0.35\hsize]{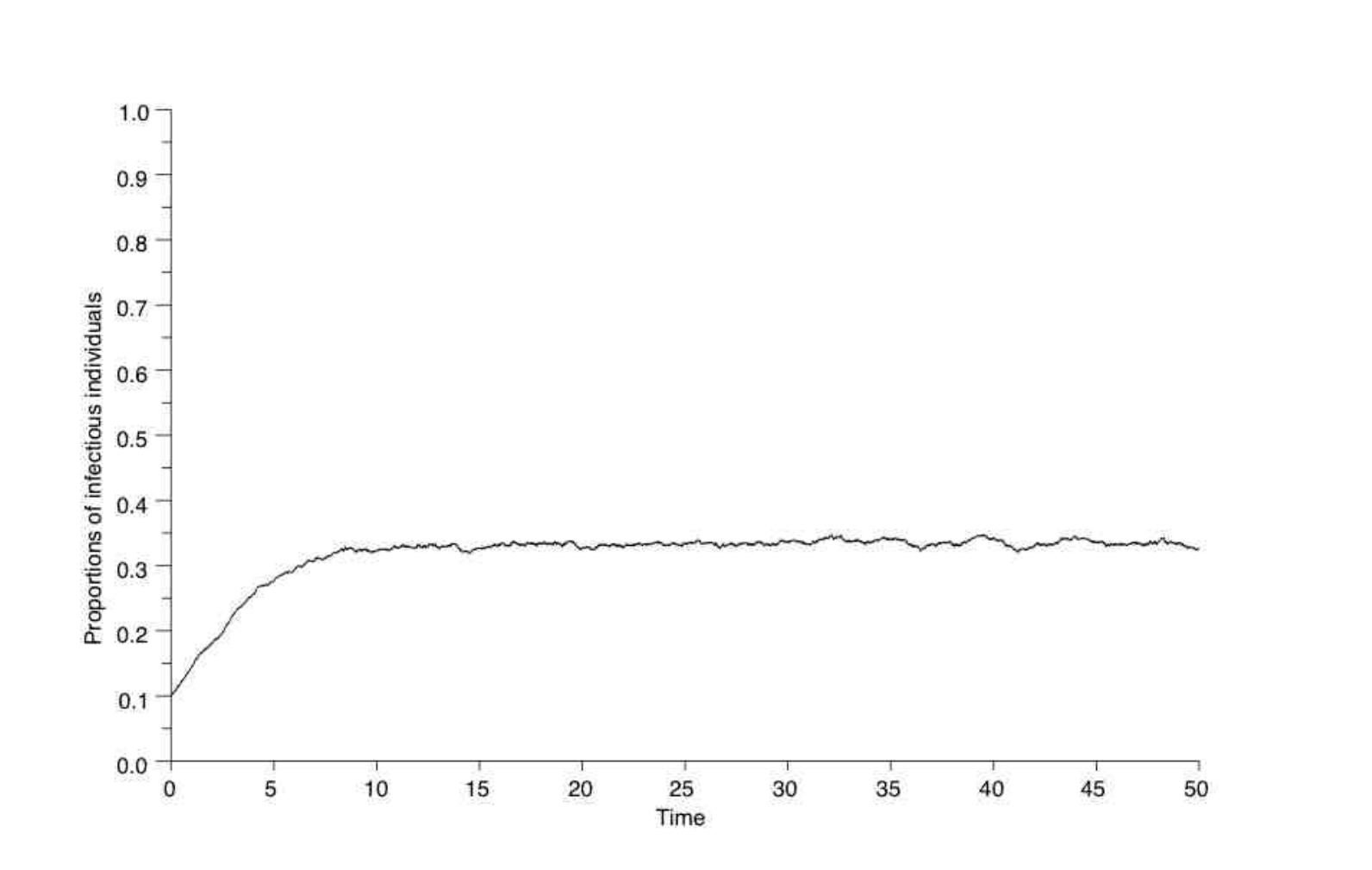}
\caption{Realized trajectories of the modified $\tau$-leaping procedure for $N=2000$ (left picture) and $N=20000$ (right picture). The parameters are as in Figure~6.}
\end{figure}

\subsection{The model with vaccination} \label{SubSecSimSIV}

We now apply the same algorithms to the model with vaccination (Example~\ref{ExampleSIV}). We consider the bistable situation from section~\ref{SubSecODESIV} (cf.~Figure~3 for the parameters) and observe the same qualitative properties concerning computation times as for the SIS model in section~\ref{SubSecSimSIS}. In Figures~8 and~9, we illustrate the realized trajectories for the SSA and the modified $\tau$-leaping algorithm, respectively (left pictures: $N=2000$, right pictures: $N=20000$). For simplicity, we only display the proportions of infectious individuals.
As in section~\ref{SubSecSimSIS}, we observe that the trajectories approach the solution of the ODE as $N \rightarrow \infty$.

\begin{figure} 
\centering \includegraphics[width=0.45\hsize,height=0.35\hsize]{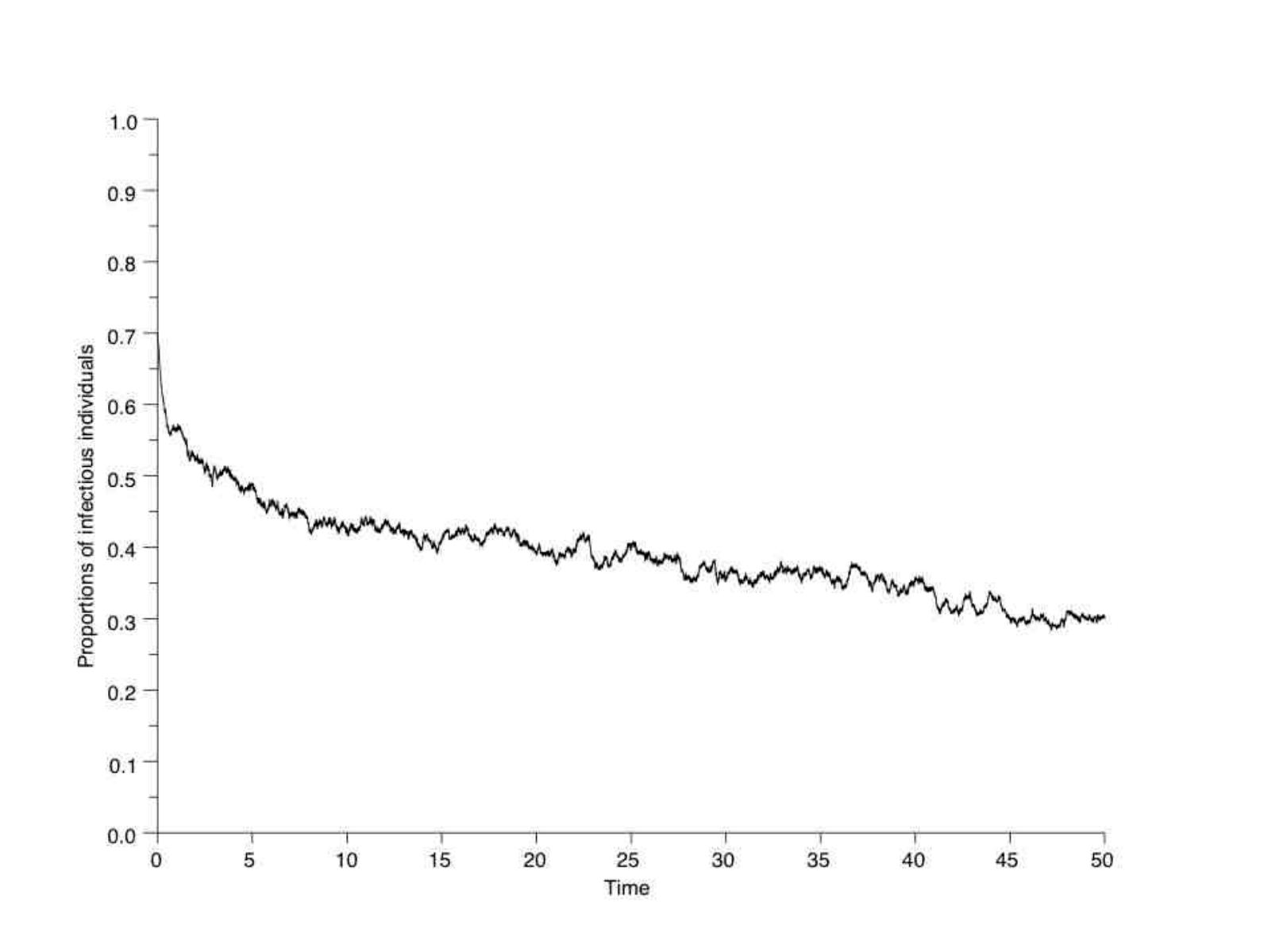}
\qquad
\includegraphics[width=0.45\hsize,height=0.35\hsize]{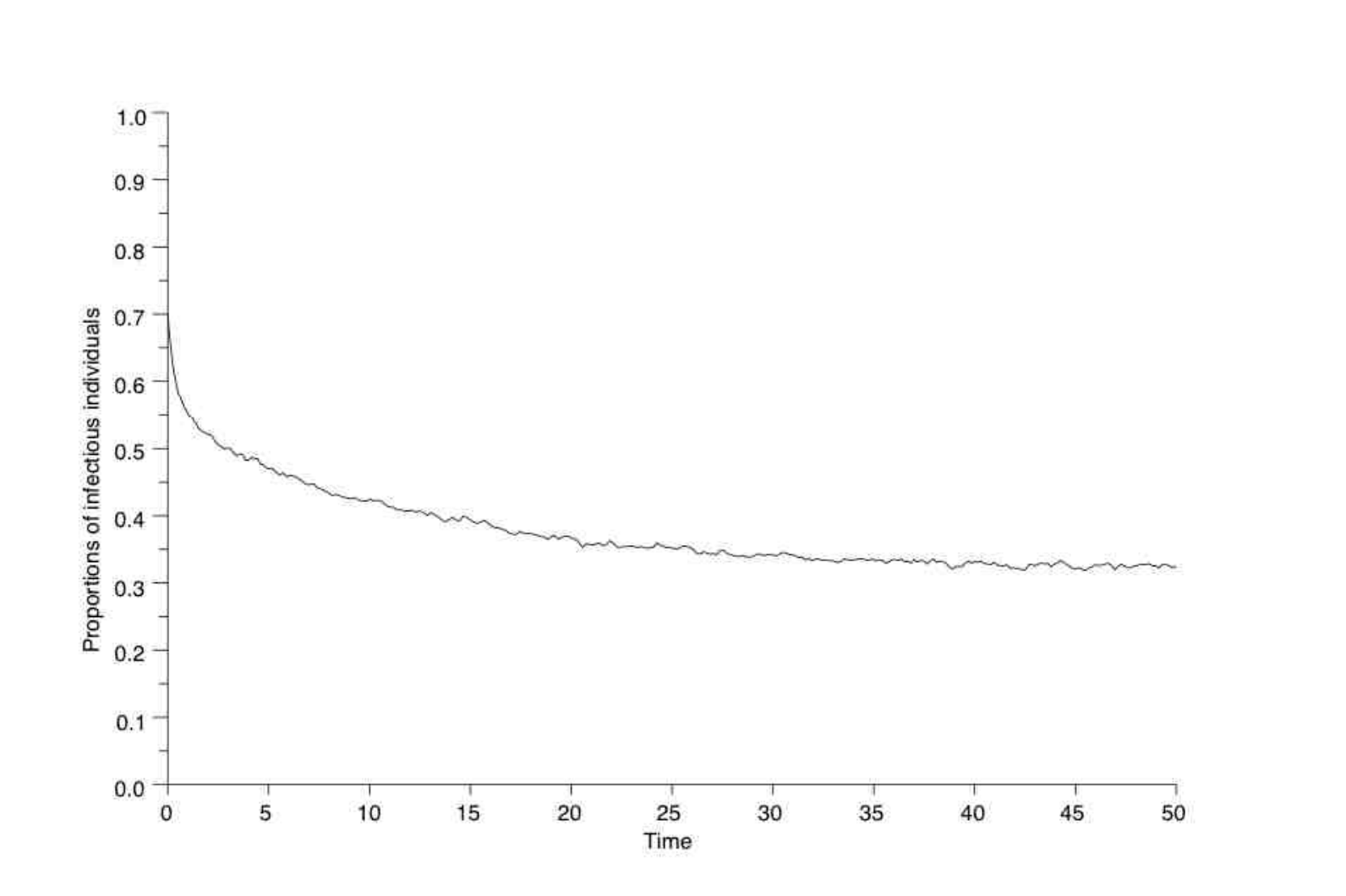}
\caption{Realized trajectories for the proportion of infectious individuals simulated by the SSA (left picture: $N=2000$, right picture: $N=20000$). $Z_2(0)=0.2$, $Z_3(0)=0.7$, $T=50$. All other parameters are as in Figure~3.}
\end{figure}

\begin{figure} 
\centering \includegraphics[width=0.45\hsize,height=0.35\hsize]{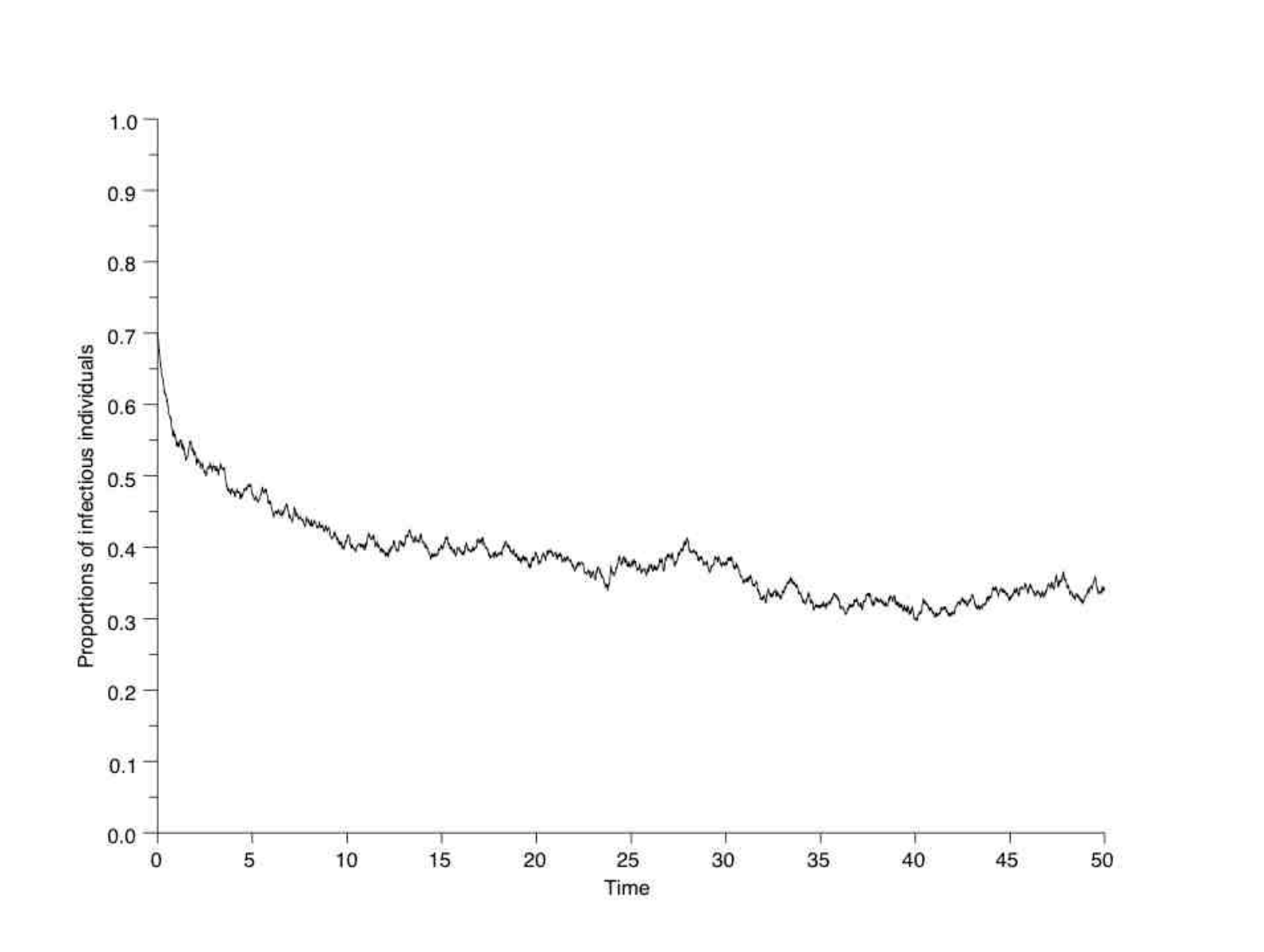}
\qquad
\includegraphics[width=0.45\hsize,height=0.35\hsize]{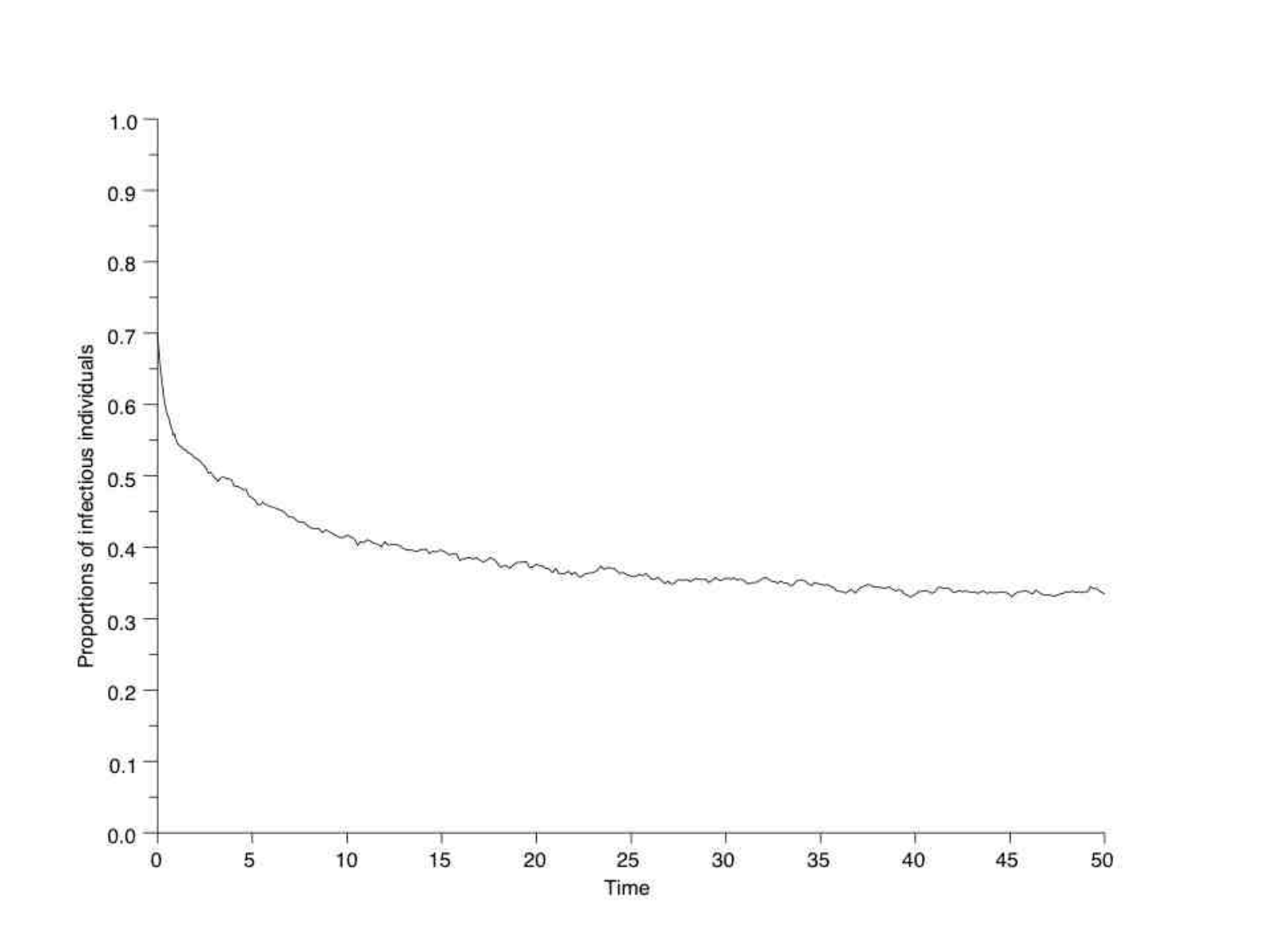}
\caption{Realized trajectories for the proportion of infectious individuals simulated by the $\tau$-leaping algorithm (left picture: $N=2000$, right picture: $N=20000$). All parameters are as in Figure~8.}
\end{figure}

It is tempting to compare those figures, in order to obtain the error induced by the tau--leaping scheme, compared to the exact SSA scheme. However, this would require that the two schemes are computed using the same realization of the Poisson processes. This is not what we have done here. We intend to do this in a near future, in order to confront simulations with the predictions of the theory concerning the error induced by the tau--leaping algorithm, see 
e.g. \cite{YHu2009}, \cite{Anderson2011}.

\subsection{Computation time}

In the next table we present the computation times (in seconds) of the SSA and modified tau-leaping for different population sizes in the SIS model. 
All programs were written and executed in Matlab. We remark that from $N=20000$
on, the computation time of the tau-leaping scheme is spectacularly smaller than that of the SSA scheme. We also note that the computation time of the tau--leaping scheme seems to be a decreasing function of $N$. The reason is probably that the larger $N$, the less frequently the scheme uses SSA steps.

\begin{center}
\begin{tabular}{|c|c|c|c|}
  \hline
  SSA & Population size $N$ & Time horizon $T$ & Computation time\\
  \hline
   & 2000 & 50 & 2.8 \\
  \hline
   & 20000 & 50 & 24.66 \\
  \hline
   & 200000 & 50 & 225.44 \\
  \hline
   Modified tau-leaping &  &  &  \\
  \hline
    & 2000 & 50 & 2.3 \\
  \hline
    & 20000 & 50 & 3.02 \\
  \hline
    & 200000 & 50 & 1.05 \\
  \hline
\end{tabular}
\end{center}

The next table presents the computation times of the SSA and modified tau-leaping for different population sizes in the SIV model. The same comments 
as for the SIS model apply.

\begin{center}
\begin{tabular}{|c|c|c|c|}
  \hline
  SSA & Population size $N$ & Time horizon $T$ & Computation time\\
  \hline
   & 2000 & 50 & 4.51 \\
  \hline
   & 20000 & 50 & 36.4 \\
  \hline
   & 200000 & 50 & 352.64 \\
  \hline
   Modified tau-leaping &  &  &  \\
  \hline
    & 2000 & 50 & 2.95 \\
  \hline
    & 20000 & 50 & 1.04 \\
  \hline
    & 200000 & 50 & 0.67 \\
  \hline
\end{tabular}
\end{center}

\subsection{Open questions} \label{SubSecSimOpen}

We conclude by describing two further open questions concerning the simulation of the process $Z^N$. We discuss briefly the simulation of models where some jump types are significantly more frequent than others (section~\ref{SubSubSecSimMalaria}). Finally, we discuss the simulation of rare events with the help of large deviations theory  (section~\ref{SubSubSecSimRare}).

\subsubsection{Models with different time scales} \label{SubSubSecSimMalaria}

The simulation algorithms described in sections~\ref{SubSecSimSSA}-\ref{SubSecSimTauMod} tend to become inefficient or inaccurate if some of the jump rates are much larger than others (cf.~also Footnote~\ref{FootnoteStiff}). In epidemiology, this is typically the case for so-called host-vector diseases (e.g., malaria), where the disease is transmitted from vectors (mosquitoes) to hosts (humans) and vice versa. In the particular case of malaria, the jumps concerning only the mosquito population have a much shorter time scale than those also concerning the human population (see, e.g., the model in~\cite{Chitnis2006}). It is hence time consuming to simulate a significant number of jumps within the human population.

In order to overcome this problem, we simulate only the jumps with small rates via $\tau$-leaping. Within the interval $[t,t+\tau[$\footnote{Here, $\tau$ is the time step chosen by the $\tau$-selection procedure applied to only the small rates.}, the jumps with large rates are approximated by the solution of the ODE (cf.~section~\ref{SecODE}).\footnote{The initial value $z$ for this computation can either be chosen explicitly, i.e., $z=Z^N(t)$ or implicitly/at the mid-point of the interval (cf.~the somehow related issue in section~\ref{SubSubSecSimTauImp}).} We intend to analyze the efficiency and accuracy of this procedure. In particular, we are interested in whether the large deviations behavior (cf.~section~\ref{SubSecLDP}) is changed by this partial ODE approximation significantly.

\subsubsection{Simulation of rare events} \label{SubSubSecSimRare}

We are particularly interested in sampling the probability of the event that the process $Z^N$ exits the domain of attraction of a stable equilibrium before a given time $T$ (cf.~section~\ref{SubSecLDP}). If $N$ is large and $T$ is relatively small, the probability of this event is typically extremely small. Applying a na\"ive Monte Carlo (MC) approach is therefore inefficient and we rely on rare events simulation techniques. In~\cite{wainrib2012}, the application of Importance Sampling is explained in the context of diffusion processes. We intend to adapt this method in the following way. Via a change of measure, we transfer the process $Z^N$ to a jump process $\tilde Z$ (see, e.g.,~\cite{Shwartz1995} for a Girsanov theorem for processes of the type of $Z^N$) for which the event of interest is much more probable. Then, we apply MC simulation to the process $\tilde Z$ and associate a weight with each realization (hence the average in the MC is replaced by a weighted average). In order to choose the process $\tilde Z$ appropriately (i.e., such that the probability of exit is ``maximized''), we use the theory of large deviations: the ``cheapest'' way to exit the domain is given by the trajectory $\phi^*$ (cf.~sections~\ref{SubSecLDP} and~\ref{SecVBar}); hence, we change the rates of the process in such a way that the corresponding rate function $I(\phi^*)$ becomes small (ideally, $\phi^*$ is approximately the solution of the ODE for the modified rates).

\begin{acknowledgement}
\textbf{Acknowledgements:}
This research was supported by the ANR project MANEGE, the DAAD, the Labex Archim\`ede and CEMRACS.

\noindent We thank Nicolas Champagnat and Tony Leli\`evre for useful discussions concerning the $\tau$--leaping algorithm.

\end{acknowledgement}

\end{document}